\documentclass[10pt]{article}
\usepackage[pctex32]{graphics}
\usepackage{amsfonts, amsmath}

\renewcommand{\O}{\mathbf{0}}

\newcommand{\vz}{\mathbf{z}}

\newcommand{\sign}{\hbox{\rm sign}}
\newcommand{\vj}{\mathbf{j}}
\newcommand{\vk}{\mathbf{k}}

\def\1{\hbox{1\kern-.15em\vrule height 1.6ex width .3pt%
\vrule width .8pt height .25pt\kern.15em}}

\def\@begintheorem#1#2{\par\bgroup{\sc #1\ #2. }\it\ignorespaces}
\def\@opargbegintheorem#1#2#3{\par\bgroup{\sc #1\ #2\ (#3). }\it\ignorespaces}
\def\@endtheorem{\egroup}

\def\endproof{{\ \vbox{\hrule\hbox{%
   \vrule height1.3ex\hskip0.8ex\vrule}\hrule
  }}\par}
\newtheorem{theorem}{{\sc Theorem}}[section]
\newtheorem{lemma}[theorem]{{\sc Lemma}}
\newtheorem{corollary}[theorem]{{\sc Corollary}}

\newtheorem{remark}[theorem]{{\sc Remark\/}}
\newtheorem{example}[theorem]{{\sc Example\/}}

\def\theequation{\arabic{equation}}
\def\@eqnnum{\hbox to .01pt{}\rlap{\rm \hskip -\displaywidth(\theequation)}}

\newcounter{rmnum}
\newenvironment{romannum}{\begin{list}{{\rm ({\it \roman{rmnum}})}}{\usecounter{rmnum}\setlength{\leftmargin}{10pt}\setlength{\itemindent}{8pt}\itemsep 0.in \parsep 0.in}}{\end{list}}

\newcounter{muni}

\newenvironment{nitemizetm}{\begin{list}{
{$\bullet$}}{\setlength{\leftmargin}{6pt}\setlength{\itemindent}{6pt}\itemsep 0.in \parsep 0.in}}{\end{list}}

\def\mojabs#1#2{\vspace{.05in}\footnotesize \parindent .2in
{\bf #1. }\ignorespaces #2\par\vspace{.1in}}

\newcommand{\Pb}{\mathbb {P}}
\newcommand{\E}{\mathbb {E}}

\newcommand{\R}{\mathbb {R}}
\newcommand{\Q}{\mathbb {Q}}
\newcommand{\Z}{\mathbb {Z}}
\newcommand{\N}{\mathbb {N}}

\title{Estimating a class of diffusions from discrete observations via approximate maximum likelihood method\thanks{This work has been partially supported by Croatian Science Foundation under the project 3526, and by Ministry of Science, Education and Sports, Republic of Croatia, Grants 037-0372790-2800 and
037058.}}

\author{Miljenko Huzak\footnote{Department of Mathematics, Faculty of Science, University of Zagreb,
Bijeni\v cka 30, HR-10002 Zagreb, Croatia ({\tt huzak@math.hr})}
}

\date{}

\begin{document}

\maketitle

\mojabs{Abstract}{
An approximate maximum likelihood method of estimation of diffusion parameters $(\vartheta,\sigma)$
based on discrete observations of a diffusion $X$ along fixed
time-interval $[0,T]$ and Euler approximation of integrals is analyzed. We assume
that $X$ satisfies a
SDE of form
$dX_t =\mu (X_t ,\vartheta )\, dt+\sqrt{\sigma} b(X_t )\, dW_t$,
with non-random initial condition. SDE is nonlinear in $\vartheta$ generally.
Based on assumption that maximum likelihood estimator
$\hat{\vartheta}_T$ of the drift parameter based on continuous observation
of a path over $[0,T]$ exists we prove that measurable estimator $(\hat{\vartheta}_{n,T},\hat{\sigma}_{n,T})$ of the parameters
obtained from discrete observations of $X$
along $[0,T]$ 
by maximization of the approximate log-likelihood function
exists, $\hat{\sigma}_{n,T}$ being consistent and asymptotically normal, and $\hat{\vartheta}_{n,T}-\hat{\vartheta}_T$ tends to zero with rate $\sqrt{\delta}_{n,T}$ in probability
when $\delta_{n,T} =\max_{0\leq i<n}(t_{i+1}-t_i )$ tends to zero with $T$ fixed. The same holds in case of an ergodic diffusion when $T$ goes to infinity in a way that $T\delta_n$ goes to zero with equidistant sampling, and we applied these to show consistency and asymptotical normality of $\hat{\vartheta}_{n,T}$, $\hat{\sigma}_{n,T}$ and asymptotic efficiency of $\hat{\vartheta}_{n,T}$ in this case.
}

\mojabs{Key words}{parameter estimation, diffusion processes, discrete observation}


\pagestyle{myheadings}
\markright{ESTIMATING A CLASS OF DIFFUSIONS VIA APPROXIMATE ML METHOD}

\normalsize

\section{Introduction}

Let $X=(X_t ,t\geq 0)$ be an one-dimensional diffusion  which satisfies
It\^{o}'s stochastic differential equation (SDE) of the 
form
\begin{equation}\begin{array}{l}
X_t =x_0 +\int_0^t\mu (X_s ,\vartheta )\, ds+
\int_0^t \sqrt{\sigma}\, b(X_s )\, dW_s ,\quad t>0. 
\end{array}\label{e1}
\end{equation}
Here, $W=(W_t ,t\geq 0)$ is an one-dimensional standard Brownian
motion,
$\mu$ and $b$ are real functions such that they ensure the 
uniqueness in law of a solution to (\ref{e1}) 
and $x_0$ is a given deterministic initial value of $X$ (see e.g.\
\cite{Yor} as a reference for SDE).

The problem is to estimate unknown vector parameter $\theta =(\vartheta , \sigma )$
of $X$, given a discrete observation $(X_{t_i},0\leq i\leq n)$ of a
trajectory $(X_t ,t\in [0,T])$ over a time interval subdivision $0=:t_0<t_1 <\cdots < t_n :=T$, ($n$
 is a positive integer) with diameter $\delta_{n,T}:=\max_{0\leq i<n}(t_{i+1}-t_i )$, $T>0$ being fixed.
Component $\vartheta$ of $\theta$ is a (vector) drift parameter, and $\sigma$ is a
diffusion coefficient parameter. We assume that $\vartheta$ belongs to
drift parameter space $\Theta$, which is an open and convex
set in Euclidean space $\R^d$, and that $\sigma$ is a positive
real number. Hence, $\theta =(\vartheta ,\sigma )$ is an element of
open and convex parameter space $\Psi$ $:=\Theta\times\langle 0,+\infty\rangle$.

Diffusion parameter estimation problems based on discrete
observations have been discussed by many authors (see
\cite{AitSahalia0,AitSahalia,AitSahalia2,Bibby,DCFZ,Dohnal,FZ,Kessler,Kloeden3,Breton,Li,Yoshida}).
Although the maximum likelihood estimator (MLE) has the usual good
properties (see \cite{DCFZ}), it may not be possible to calculate it
explicitly because the transition density of process $X$ is
generally unknown and so the likelihood function (LF) of the
discrete process is unknown as well. Hence, other methods of
estimations have to be considered.

The method of parameter estimation which is discussed in this paper
and described in Section 3 below,
is based on a Gaussian approximation of the transition density and can be interpreted as based
on maximization of a discretized continuous-time log-likelihood function (LLF) as well.
Such methods
are usually called quasi-likelihood or approximate maximum likelihood (AML) methods, and in these ways
obtained estimators we will briefly call approximate maximum likelihood estimators (AMLEs).

Motivation for analyzing the method described in Section 3 is in the fact that it can provide us with
useful estimators of the parameters. It is well known that in a such way obtained AMLE of diffusion coefficient parameter $\sigma$ is consistent and asymptotically normally distributed over fixed observational time interval $[0,T]$ when $\delta_{n,T}\rightarrow 0$  (see \cite{Dohnal} in case where all drift parameters are known, and see \cite{Jacod} in general cases). The same holds in ergodic
diffusion cases when $T\rightarrow +\infty$ in a way that $\delta_{n,T}=T/n\rightarrow 0$
for appropriate equidistant sampling (see e.g.\ \cite{FZ} or \cite{Kessler}).
Local asymptotic properties of the AMLE of drift parameters over fixed interval $[0,T]$ and when $\delta_{n,T}\rightarrow 0$ are less known especially in more general cases,
particularly when drift is nonlinear in its parameters (see \cite{Bishwal}).
Although a knowledge of local asymptotic properties of drift parameter AMLEs does not imply their consistency or asymptotic normality necessarily it may help in further analysis of the AMLEs which might include, for example, measuring effects of discretization on the estimator's standard errors with applications in simulation studies. In ergodic diffusion cases it is well known that the AMLE of drift (vector) parameter is consistent and asymptotically normal and efficient when $T \rightarrow +\infty$ in a way that $T\delta_{n,T}^2\rightarrow 0$ for equidistant sampling (see e.g.\ \cite{FZ} for one-dimensional case and \cite{Kessler} for vector and more general cases) but the rate of convergence
of $\hat{\vartheta}_{n,T}-\hat{\vartheta}_T$ to zero are still less investigated.  Let us stress that the
problems of statistical inferences about diffusion drift parameters are very important especially in biomedical modeling
(see \cite{HuzakI}).

For the completeness we should also stress that local convergence of the AMLE of both vector
parameters $\theta =(\vartheta,\sigma )$ to the MLE of $\theta$ based on discrete observations and equidistant sampling,
have been investigated (see \cite{AitSahalia0, AitSahalia2,Li}). Let $\tilde{\theta}_{n,\delta}$ denote MLE of
$\theta$ based on discrete observations with $\delta_{n,T}\equiv\delta =$ const., and  let $\tilde{\theta}_{n,\delta}^{(k)}$ be AMLE obtained from an approximate LF based on a closed-form $k$th order approximation of the transition densities.
Then in case of Hermite-polynomial-based analytical expansion approach for approximation of transition density, $\tilde{\theta}_{n,\delta}^{(k_n)}-\tilde{\theta}_{n,\delta}\rightarrow 0$ when $k_n\rightarrow\infty$, and a sequence $(k_n)$
can be chosen sufficiently large to deliver any rate of convergence (see \cite{AitSahalia0}),
and there exist sequences of regular matrices $(S_{n,\delta} )$ and positive numbers $(\delta_n)$ such that $\delta_n\rightarrow 0$ and
$S_{n,\delta_n}^{-1}(\tilde{\theta}_{n,\delta_n}^{(k)}-\tilde{\theta}_{n,\delta_n})=O_{\Pb}(1)$ (see \cite{AitSahalia2}).
For an alternative approach to approximation and analog results, see \cite {Li}.

In this paper we analyze the considered AMLE of drift parameters by studying the relation between the AMLE and the
MLE obtained from continuously observed diffusion paths. We state general conditions for proving and prove: (1.) existence and measurability of the AMLE, (2.)
that $\hat{\vartheta}_{n,T}-\hat{\vartheta}_T$ converges to zero with rate
$\sqrt{\delta_{n,T}}$ in probability
when $\delta_{n,T}\rightarrow 0$ over fixed bounded observational time interval $[0,T]$,
and (3.) that $\hat{\vartheta}_{n,T}-\hat{\vartheta}_T$ converges to zero with rate
$\sqrt{\delta_{n,T}}$ in probability
when $T\rightarrow +\infty$ in a way that $\delta_{n,T}=T/n \rightarrow 0$ in an ergodic diffusion case and equidistant sampling. We apply these findings in
proving: (4.) measurability, consistency and asymptotic normality of diffusion coefficient parameter AMLEs when $\delta_{n,T}\rightarrow 0$ in both cases: when $T$ is fixed, and in an ergodic diffusion case when $T\rightarrow +\infty$ and $T\delta_{n,T}=T^2/n\rightarrow 0$ with equidistant sampling, and (5.)
consistency and asymptotic normality and efficiency of drift parameter AMLEs in an ergodic case when $T\rightarrow +\infty$ in a way that $T\delta_{n,T}\rightarrow 0$ with equidistant sampling.

Properties (1.-2.) for drift parameter AMLEs were proved
in \cite{Breton} in cases when drift depended linearly on its parameters. For detailed review of liner case see \cite{Bishwal}. The first nonlinear case was covered by the author in his Ph.D. thesis \cite{Huzak}. The main assumption was that the drift was an analytic function in its parameters with properly bounded derivatives of all orders. In this paper we only assume that the drift has at least $d+3$ continuous derivatives with respect to the drift parameters ($d$ is a dimension of the drift parameter vector). The main difficulty was in proving core technical Theorem \ref{tm:new} of Section 5. Although facts (4.-5.) have been already known we included these alternative proofs  for completeness and the illustrative purposes of the applicability of the findings (1.-3.) and in this paper developed methods. We belive that other discretization schemes (for example, of higher order) can be analyzed similarly by using the techniques of this paper.

The paper is organized in the following way. In the next section we introduce
notation used through the paper. The discussed method of
estimation is described in Section 3. The main results are presented in Section 4.
Examples are provided in Section 5.  
The proofs of the main results are  in the last section. Lemmas are proved in Appendix.

\section{Notations}

Let
$|\cdot |$ denote Euclidean norm in $\R^d$ and its induced operator norm, and let
$|\cdot |_\infty$ be max-norm. If $f$ is a bounded
real function, $\|f\|_\infty:=\sup_\vz |f(\vz)|$ is a sup-norm of $f$. Let $L^p (\Pb )$ be
the Banach space of all random variables with finite $p$-th moment 
and let $\|\cdot\|_{L^p(\Pb )}$ denote its 
norm.

If $(x,\vartheta )\mapsto f(x,\vartheta )$ 
is a real function defined on an open subset of $\R\times\R^d$,
then we denote by $D_{\vartheta }^m f(x,\vartheta )$ the $m$-th partial derivative
with respect to $\vartheta$. Let
$|D_{\vartheta }^m f(x,\vartheta ) |_\infty :=
\max_{j_1 +\cdots +j_d =m}|\frac{\partial^m f}{\partial\vartheta_1^{j_1} \ldots\partial
\vartheta_d^{j_d}}|$. In this case we say that $D_{\vartheta }^m f(x,\vartheta )$ is bounded
if all partial derivatives $\frac{\partial^m f}{\partial\vartheta_1^{j_1} \ldots\partial
\vartheta_d^{j_d}}(x,\vartheta )$ 
are bounded, and $\|D_{\vartheta }^m f\|_\infty
:=\max_{j_1 +\cdots +j_d =m}\|\frac{\partial^m f}{\partial\vartheta_1^{j_1} \ldots\partial
\vartheta_d^{j_d}}\|_\infty$.
The notation $D^2_{\vartheta } f(x,\vartheta )<\O$ means that the
Hessian $D^2_{\vartheta } f(x,\vartheta )$ is a negatively definite matrix. Similarly for a positively definite matrix.
$D_{\vz}^0 f\equiv f$ by convention.
The $m$-th derivative of $f$ at a point $\vz$ we simply denote by $D^m f(\vz)$.

Let ${\cal K}$ and $\Theta$ be open sets in $\R^d$.
The closure and the boundary of ${\cal K}$ will be denoted by
$\overline{\cal K}$ and $\partial {\cal K}$ respectively, and the $\sigma$-algebra of Borel subsets of $\Theta$ by
${\cal B}(\Theta )$. If ${\cal K}\subset\Theta$ is an open set such that $\overline{\cal K}$ is compact in $\Theta$
then we will say that ${\cal K}$ is a relatively compact set in $\Theta$.

Let $(\gamma_n ,n\geq 1 )$ be a sequence of positive numbers and let
$(Y_n ,n\geq 1 )$ be a sequence of random variables defined on some probability space.
We will say
that $(Y_n ,n\geq 1 )$ is $O_\Pb (\gamma_n)$,  and write
$Y_n =O_\Pb (\gamma_n )$, if the sequence $(Y_n /\gamma_n ,n\geq 1)$
is bounded in probability, i.e.\ if
\[\begin{array}{l}
\lim_{A\rightarrow +\infty}\overline{\lim_n}\,\Pb\{ \gamma_n^{-1}|Y_n|> A\} =0.
\end{array}\]

\section{Estimation method}

Let $0=t_0<t_1<\cdots <t_n=T$ be discrete times at which diffusion $X$
is observed, and let us denote by $\Delta$  the difference operator
defined in the following way: if $F$ is a function defined on
$[0,T]$ then $\Delta_i F:=F(t_{i+1})-F(t_i)$, $0\leq i<n$.

Let us discretize SDE (\ref{e1}) over interval $[t_i,t_{i+1}]$ by
using the Euler approximation of the both types of integrals:
\[\begin{array}{l}
X_{t_{i+1}}-X_{t_i}\approx\mu (X_{t_i},\vartheta )
(t_{i+1}-t_{i})+\sqrt{\sigma}\, b(X_{t_i}) (W_{t_{i+1}}-W_{t_i}).
\end{array}\]
In this way the following stochastic difference equation is
obtained:
\begin{equation}\begin{array}{l}
\Delta_i Z=\mu (Z_i ,\vartheta )\, \Delta_it+\sqrt{\sigma}\, b(Z_i
)\,\Delta_i W
\end{array}\label{e2}
\end{equation}
for $0\leq i<n$, and $Z_0 =x_0$.
Solution to (\ref{e2}) is a time-discrete process
$Z=(Z_0,Z_1,\ldots,Z_n)$ that is an approximation of $X$ over
$[0,T]$.
Up to the constant not depending on the parameters 
a LLF of the process $Z$ is
\begin{equation}\begin{array}{l}
-\frac{1}{2}\sum_{i=0}^{n-1}\left(\frac{
 (\Delta_i Z-\mu (Z_i ,\theta )\Delta_i t)^2}{\sigma b^2 (Z_i )\Delta_i t}+
\log\sigma\right).
\end{array}\label{LLFd}
\end{equation}
Criterion function
\begin{equation}\begin{array}{l}
{\cal L}_{n,T} (\theta )={\cal L}_{n,T} (\vartheta ,\sigma ):=
-\frac{1}{2}\sum_{i=0}^{n-1}\left(\frac{(\Delta_i X-\mu (X_{t_i} ,
\vartheta )\Delta_it )^2}{\sigma b^2 (X_{t_i} )\Delta_it}+
\log\sigma\right)
\end{array}\label{e3}
\end{equation}
is obtained from (\ref{LLFd}) by substituting $(Z_i ,0\leq i\leq n)$
with discrete observations $(X_{t_i} ,0\leq i\leq n)$ of
diffusion $X$. Notice that
\[\begin{array}{l}
{\cal L}_{n,T} (\vartheta ,\sigma )=-\frac{1}{2\sigma}\sum_{i=0}^{n-1}
\frac{(\Delta_i X)^2}{b^2 (X_{t_i}) \Delta_it}-\frac{n}{2}\log\sigma +
\frac{1}{\sigma}{\ell}_{n,T} (\vartheta) ,
\end{array}\]
where
\begin{equation}\begin{array}{l}
\ell_{n,T} (\vartheta )=\sum_{i=0}^{n-1} \frac{\mu (X_{t_i},\vartheta )}{b^2
(X_{t_i})}\Delta_i X-\frac{1}{2}\sum_{i=0}^{n-1}
\frac{\mu^2 (X_{t_i},\vartheta )}{b^2 (X_{t_i})}\Delta_it \label{lnizraz}
\end{array}\end{equation}
depends only on drift parameter $\vartheta$.

A point of maximum $\hat{\theta}_{n,T} =(\hat{\vartheta}_{n,T} ,
\hat{\sigma}_{n,T} )$ of function (\ref{e3}) in $\Psi$  
is an AMLE of vector parametar $\theta$ if it exists. Notice that if AMLE
exists then necessary
\begin{equation}
D{\cal L}_{n,T} (\hat{\vartheta}_{n,T},\hat{\sigma}_{n,T} )=0
\,\Leftrightarrow\,\left\{
\begin{array}{l}  D\ell_{n,T} (\hat{\vartheta}_{n,T})=0\\
\hat{\sigma}_{n,T} =\frac{1}{n}\sum_{i=0}^{n-1}\frac{(\Delta_i X-\mu (X_{t_i},\hat{\vartheta}_{n,T})\Delta_it )^2
}{b^2 (X_{t_i})\Delta_it}.
\end{array}\right. \label{defv}
\end{equation}
Hence every stationary point $\hat{\vartheta}_{n,T}$ of function $\ell_{n,T}$
uniquely determines second component $\hat{\sigma}_{n,T}$ of stationary point $\hat{\theta}_{n,T}=(\hat{\vartheta}_{n,T},\hat{\sigma}_{n,T})$ of
function ${\cal L}_{n,T}$
by the following expression:
\begin{equation}\begin{array}{l}
\hat{\sigma}_{n,T} =\frac{1}{n}\sum_{i=0}^{n-1}\frac{(\Delta_i X-\mu (X_{t_i},
\hat{\vartheta}_{n,T} )\Delta_it )^2}{b^2 (X_{t_i})\Delta_it}.
\end{array}\label{hatvn}
\end{equation}
Moreover,
if $\hat{\vartheta}_{n,T}$ is a unique point of the global maximum of
function $\ell_{n,T}$ then $\hat{\theta}_{n,T}$ is a unique point of the global
maximum of function ${\cal L}_{n,T}$. Hence to prove existence of a measurable AMLE $\hat{\theta}_{n,T}$
it is sufficient to prove that there exists a measurable point of maximum of function $\ell_{n,T}$.

\section{Main results}

\subsection{{\em Fixed maximal observational time case\/}}

Let the following assumptions be satisfied.
\vspace{5pt}

{\sc (H1a):} For all $\theta =(\vartheta ,\sigma )\in\Psi$,
there exists a strong solution $(X,W)$ of the SDE (\ref{e1}) on time interval
$[0,+\infty\rangle$
with values in open interval
$E\subseteq\R$. \vspace{5pt}

{\sc (H2a):} For all $\vartheta\in\Theta$, $\mu (\cdot,\vartheta )\in C^2(E)$ 
and 
$b\in C^3 (E)$. 
Moreover for all $x\in E$, $b(x)\neq 0$ and $\sign\, b =$ const.
\vspace{5pt}

For example, by Theorem 5.2.2 in \cite{Friedman}, 
{\sc (H1a)} will be satisfied if in addition to 
{\sc (H2a)} we assume that for all $\vartheta\in\Theta$  SDE (\ref{e1})
satisfies so called the bounded linear growth assumption, i.e.\ that
there exists a positive constant $C$ such that for all $x\in E$,
$ 
|\mu (x,\vartheta )|+|b(x)|\leq C(1+|x|).
$ 
More precisely, {\sc (H2a)} states that the functions $x\mapsto b(x)$ and
$x\mapsto\mu (x,\vartheta )$, $\vartheta\in\Theta$, are continuously
differentiable in $E$ and hence locally Lipschitz. In this case there
exists a strong, continuous and pathwise unique solution to SDE (\ref{e1})
on time interval $[0,+\infty\rangle$.
However, there are some SDEs which satisfy {\sc (H1a)} and {\sc (H2a)} but do not
satisfy the linear growth assumption
(see e.g. Example \ref{primjer1} of Section 5). 
\vspace{5pt}

{\sc (H3a):} For all $(x,\vartheta )\in E\times\Theta$ and all $1\leq m\leq d+3$, there
exists partial derivatives $D_\vartheta^m\mu (x,\vartheta )$,
$\frac{\partial}{\partial x}D_\vartheta^m\mu (x,\vartheta )$,and
$\frac{\partial^2}{\partial x^2}D_\vartheta^m\mu (x,\vartheta )$
of drift function $\mu$.
Moreover, for all $0\leq m\leq d+3$,
$D_\vartheta^m\mu$, $\frac{\partial}{\partial x}D_\vartheta^m\mu$, $\frac{\partial^2}{\partial x^2}D_\vartheta^m\mu\in C(E\times\Theta)$.
\vspace{5pt}


Let $\Pb_\theta$ denote the law of $X$ for $\theta\in\Psi$.
We assume that probabilities $\Pb_\theta$, $\theta
\in\Psi$, are defined on filtered space
$(\Omega ,({\cal F}_T^0, T\geq 0) )$ where $\Omega$ is a set of continuous functions $\omega :[0,+\infty\rangle\rightarrow
E$ such that $\omega (0)=x_0$, ${\cal F}_T^0$ is a
$\sigma$-algebra generated by the coordinate functions up to the time $T$, and the filtration is augmented in so called the usual way
(see e.g.\ I.4 in \cite{Yor}).
On this space, coordinate process $(\omega\mapsto\omega (t),
t\geq 0)$ is a canonical version of $X$ (see \cite{Yor}, I.\S3).
Hence, for each $T>0$ we assume that $X$ is defined on the measurable space $(\Omega ,
{\cal F}_T^0 )$ as a canonical process with  law $\Pb_\theta$.

For the moment, let us assume that we are able to observe the process
$(X_t ,0\leq t\leq T)$ continuously. Because diffusion coefficient
parameter $\sigma$ can be uniquely determined through equation
\begin{equation}\begin{array}{l}
\sigma =\frac{\lim_n \sum_{i=1}^{2^n}(X_{jT2^{-n}}-X_{(j-1)T2^{-n}})^2}{
\int_{0}^{T}b^2 (X_t )\, dt}\quad (\mbox{\rm a.s.}\;\Pb_\theta )
\end{array}\label{sigma}
\end{equation}
(see \cite{BH}) since $b^2>0$ by {\sc (H2a)}, the estimation problem from continuously
observed process can be reduced to an estimation problem for drift
parameter $\vartheta\in\Theta$. In this case
for every fixed diffusion parameter $\sigma$ assumed to be known,
and every two different $\vartheta_1 ,\vartheta_2\in\Theta$, probability
measures $\Pb_{(\vartheta_1 ,\sigma )}$ and $\Pb_{(\vartheta_2 ,\sigma )}$
are equivalent on ${\cal F}_T^0$, and
\[\begin{array}{l}
\log\frac{d\Pb_{(\vartheta_2 ,\sigma )}}{d\Pb_{(\vartheta_1 ,\sigma )}}\! =\!
\frac{1}{\sigma}(\!\int_0^T\!\!\!\frac{\mu (X_t ,\vartheta_2 )\! -\!
\mu (X_t ,\vartheta_1 )}{b^2 (X_t )} dX_t\! -\!\frac{1}{2}
\int_0^T\!\!\!\frac{\mu^2 (X_t ,\vartheta_2 )\! -\!\mu^2 (X_t ,\vartheta_1 )}{
b^2 (X_t )} dt)
\end{array}\]
where
$\frac{d\Pb_{(\vartheta_2 ,\sigma )}}{d\Pb_{(\vartheta_1 ,\sigma )}}$
denotes Radon-Nikodym derivative of $\Pb_{(\vartheta_2 ,\sigma)}$
with respect to $\Pb_{(\vartheta_1 ,\sigma)}$ on ${\cal F}_T^0$ (see \cite{Feigin}).
If we fix some $\vartheta_*\in\Theta$, a continuous-time LLF is
$\vartheta\mapsto\log\frac{d\Pb_{(\vartheta ,\sigma)}}{
d\Pb_{(\vartheta_* ,\sigma )}}$. Up to the constant and factor not depending on $\vartheta$,
function
\begin{equation}\begin{array}{l}
\ell_T (\vartheta ):=
\int_0^T\frac{\mu (X_t ,\vartheta )}{b^2 (X_t )}\, dX_t -
\frac{1}{2}\int_0^T\frac{\mu^2 (X_t ,\vartheta )}{b^2 (X_t )}\, dt.
\end{array}\label{ell}
\end{equation}
is equal to the LLF.
Hence, $\ell_T$ will be called a
continuous-time LLF (see \cite{Lanska}). Assumption
{\sc (H3a)} implies that $\ell_T$ is at least three-times continuously
differentiable function on $\Theta$, and for $1\leq m\leq d+3$, its  derivatives are equal
to (see \cite{Lanska} for $m\leq 2$)
\begin{equation}
\begin{array}{l}
D^m\ell_T (\vartheta )=
\int_0^T\frac{1}{b^2 (X_t )}D_\vartheta^m\mu (X_t ,\vartheta )\, dX_t -
\frac{1}{2}\int_0^T\frac{1}{b^2 (X_t )}D_\vartheta^m\mu^2 (X_t ,\vartheta )\,
dt.
\end{array}\label{Dell}
\end{equation}

{\sc (H4a):} For all $\omega\in\Omega$, function
$\vartheta\mapsto\ell_T (\vartheta )=\ell_T (\vartheta, \omega )$
has a unique point of global maximum
$\hat{\vartheta}_T=\hat{\vartheta}_T(\omega )$ in $\Theta$.
 Moreover, $D_{\vartheta}^2 \ell_T (\hat{\vartheta}_T)<\O$.
\vspace{5pt}

Assumption {\sc (H4a)} 
enables property ($ii$) in
Theorem \ref{tm:30} below, to be proved. If {\sc (H3a)} and {\sc (H4a)} hold then Lemma 4.1.\
from \cite{huzakMLE} implies that $(\omega,\vartheta)\mapsto\ell_T(\vartheta)(\omega)$ is
an ${\cal F}_T^0\otimes {\cal B}(\Theta )$-measurable function, and continuous-time MLE
$\hat{\vartheta}_T$ is an ${\cal F}_T^0$-measurable random variable.

Let ${\cal F}_{n,T}$ be a $\sigma$-subalgebra of ${\cal F}_T^0$ generated by discrete observation
$(X_{t_i}, 0\leq i\leq n)$  of process
$(X_t ,0\leq t\leq T)$. Notice that if {\sc (H3a)} holds then $(\omega,\vartheta )\mapsto{\ell}_{n,T} (\vartheta, \omega )$ (given by \ref{lnizraz})
is an ${\cal F}_{n,T}\otimes {\cal B}(\Theta)$ measurable function by Lemma 4.1.\ in \cite{huzakMLE}.

If $\ell_{n,T}$ is a concave function on $\Theta$ then a stationary point $\hat{\vartheta}_{n,T}$ is an unique point of maximum of $\ell_{n,T}$ on $\Theta$ and hence it is ${\cal F}_{n,T}$-measurable  by e.g. Lemma 4.1.\ in \cite{huzakMLE}.
If $\ell_{n,T}$ is not a concave function on $\Theta$,
for proving ${\cal F}_{n,T}$-measurability of estimators
$\hat{\vartheta}_{n,T}$ (and so $\hat{\theta}_{n,T}$) introduced in Section 3 we
need additional assumptions:
\vspace{5pt}

{\sc (H5a):} $\Theta$ is a relative compact set in $\R^d$, and for each $0\leq m\leq d+3$,
$D_\vartheta^m\mu$, $\frac{\partial}{\partial x}D_\vartheta^m\mu$, $\frac{\partial^2}{\partial x^2}D_\vartheta^m\mu \in C(E\times\overline{\Theta})$.\vspace{5pt}

{\sc (H6a):}
For all $\omega\in\Omega$ and some $r>0$,
\[\begin{array}{l}
\ell_T (\hat{\vartheta}(\omega ), \omega)>\sup_{|x|\geq r}\ell_T (\hat{\vartheta}(\omega )+x, \omega).
\end{array}\]

\noindent
Assumption {\sc (H6a)} holds if {\sc (H5a)} holds and $\hat{\vartheta}_T$ is the unique point of maximum of
$\ell_T$ on compact $\overline{\Theta}$.

\begin{theorem}\label{tm:30}
Let us assume that {\sc (H1a-4a)} hold and $T>0$ be fixed.
Then there exists a sequence $(\hat{\vartheta}_{n,T}, n\geq 1 )$ of ${\cal F}_T^0$-measurable random vectors
such that for all $\theta =(\vartheta ,\sigma )\in\Psi$
and when $\delta_{n,T}\downarrow 0$,
\begin{romannum}
\item $\lim_n\Pb_{\theta}(D{\ell}_{n,T}
(\hat{\vartheta}_{n,T})=\O )=1$
\item $(\Pb_{\theta})\lim_n\hat{\vartheta}_{n,T} =
\hat{\vartheta}_T$
\item $\hat{\vartheta}_{n,T} -\hat{\vartheta}_T=O_{\Pb_\theta}(\sqrt{\delta_{n,T}})$, $n\rightarrow +\infty$
 \item If $(\tilde{\vartheta}_{n,T},n\geq 1)$ is an ${\cal F}_T^0$-measurable sequence in $\Theta$ that satisfies $(i-ii)$ then
 $\lim_n\Pb_{\theta} (\tilde{\vartheta}_{n,T} =\hat{\vartheta}_{n,T})=1$.
\end{romannum}
If either for $n\geq 1$ and
almost all $\omega\in\Omega$ function $\vartheta\mapsto
{\ell}_{n,T}(\vartheta, \omega )$ has a unique point of local
maximum which is a point of the global maximum as well,
or the hypotheses  {\sc (H5a-6a)} are satisfied,
then $\hat{\theta}_{n,T}$ can be chosen to be
${\cal F}_{n,T}$-measurable.
\end{theorem}

\begin{corollary}\label{cor:30a}
Let {\sc (H1a-4a)} hold, $T>0$ be fixed, and  $(\hat{\sigma}_{n,T},n\geq 1)$
be given by $(\ref{hatvn})$. Then
\begin{romannum}
\item $(\Pb_{\theta})\lim_n\hat{\sigma}_{n,T}
=\sigma$;
\item $(\sqrt{n}\frac{1}{\sigma\sqrt{2}
}(\hat{\sigma}_{n,T} -\sigma),n\geq 1)$ converges in law w.r.t.\ $\Pb_\theta$ to the standard normal
distribution $N(0,1)$ when $n\rightarrow +\infty$.
\end{romannum}
Moreover, if $\hat {\vartheta}_{n,T}$ is ${\cal F}_{n,T}$-measurable then $\hat{\sigma}_{n,T}$ is ${\cal F}_{n,T}$-measurable too.
\end{corollary}

\begin{remark}\label{rem1}{\rm
Theorem \ref{tm:30} still holds if we replace {\sc (H1a)} with the assumption that $T<\xi$ a.s. where $\xi$ is a maximal random time such that SDE (\ref{e1}) has a solution on $[[0,\xi[[=\{(\omega,t)\in\Omega\times [0,+\infty\rangle\,:\;0\leq t<\xi (\omega)\}$.
$\xi$ exists by assumption {\sc (H2a)} and the existence and uniqueness theorem for SDEs (see e.g.\ \cite{Friedman} or \cite{Yor}).
}
\end{remark}

\begin{remark}\label{rem2}{\rm
Theorem \ref{tm:30} still holds if the drift and diffusion coefficient functions depend on time variable too (non autonomous case: $(t,x)\mapsto\mu (t,x,\vartheta )$, $\sigma b(t,x)$) in a way that assumptions {\sc (H2a)} and {\sc (H3a)} hold for $\mu$ and $b$
with $x$ and $E$ replaced with $(t,x)$ and $\tilde{E}=[0,+\infty\rangle\times E$ respectively.
}
\end{remark}

\subsection{{\em Ergodic diffusions case}}

Let the coefficient diffusion function parameter $\sigma >0$ be fixed. We need
the following assumptions.\vspace{5pt}

{\sc (H1b):} {\sc (H1a)} holds, and $X$ is an ergodic diffusion with stationary distribution $\pi_\vartheta (dx)$,
 $\vartheta\in\Theta$. \vspace{5pt}

{\sc (H2b):} {\sc (H2a)} holds, and for all $\vartheta\in\Theta$ functions $\mu (\cdot ,\vartheta )b'/b, (b')^2, b''b\in L^{16} (\pi_\vartheta )$, $b^2b'''\in L^{8} (\pi_\vartheta )$, and there exist a function $c\in L^{1} (\pi_\vartheta )$ and a number $h_0>0$ such that 
\[\begin{array}{l}
\sup_{0<h\leq h_0}\E_{(\vartheta,\sigma)}\exp\left( 8\int_0^{h}\left(2\frac{\mu (\cdot ,\vartheta ) b'}{b} +\sigma (b''b+15b'^2)\right)(X_s)\,ds\right)\leq c(x_0).
\end{array}\]

{\sc (H3b):} {\sc (H3a)} and {\sc (H5a)} hold, and there exist nonnegative functions $g_{0},g_1,g_2:E\rightarrow\R$ such that
for all $\vartheta_0\in\Theta$, $g_0\in L^{32} (\pi_{\vartheta_0} )\cap C^1(E)$ such that $g_0' b\in L^{16} (\pi_{\vartheta_0} )$,
$g_1\in L^{16} (\pi_{\vartheta_0} )\cap C (E)$, $g_2\in L^{8} (\pi_{\vartheta_0} )\cap C (E)$, and for all $x\in E$ and $0\leq m\leq d+3$,
\[\begin{array}{lcl}
\sup_{\vartheta\in\overline{\Theta}}|D^m_\vartheta\mu (x,\cdot )/b(x)|_\infty &\leq &g_0(x)\\
\sup_{\vartheta\in\overline{\Theta}}|\frac{\partial}{\partial x}D^m_\vartheta\mu (x,\cdot )|_\infty
&\leq & g_1(x)\\
\sup_{\vartheta\in\overline{\Theta}}|\frac{\partial^2}{\partial x^2}D^m_\vartheta\mu (x,\cdot )b(x)|_\infty &\leq & g_2(x).
\end{array}\]

{\sc (H4b):} For all $\vartheta\in\Theta$,
\begin{equation}\begin{array}{l}
(\forall \vartheta'\in\overline{\Theta})\; \vartheta'\neq\vartheta\Rightarrow
\int_E \frac{(\mu (x,\vartheta)-\mu(x,\vartheta'))^2}{b^2 (x)}\, \pi_\vartheta (dx) >0.
\end{array}\label{identiassump}
\end{equation}

{\sc (H5b):} For all $\vartheta\in\Theta$, functions $\frac{\partial\mu }{\partial\vartheta_i}(\cdot,\vartheta )/b$, $1\leq i\leq d$, are linearly independent in $L^2 (\pi_\vartheta )$.
\vspace{5pt}

\noindent
$\Theta$ is a relatively compact set in $\R^d$ by assumption {\sc (H5a)} since {\sc (H3b)} holds.
Assumptions {\sc (Hb1-b3)} imply that for all $\vartheta_0\in{\Theta}$ and $\vartheta\in\overline{\Theta}$, $\Pb_{(\vartheta_0,\sigma)}$-a.s.
\begin{equation}\begin{array}{l}
\lim_{T\rightarrow +\infty}\frac{1}{T}\ell_T (\vartheta )=\frac{1}{2}\int_E\frac{\mu (x,\vartheta_0)^2 -
(\mu (x,\vartheta_0 )-\mu (x,\vartheta ))^2}{b^2(x)}\,\pi_{\vartheta_0} (dx)
=:\ell_{\vartheta_0} (\vartheta )
\end{array}\label{liml0}
\end{equation}
by ergodic property of the diffusion and the law of large numbers for continuous martingales (see e.g. \cite{Yor}, Chapters V and X).
Function $\ell_{\vartheta_0} :\overline{\Theta}\rightarrow\R$ defined for every $\vartheta_0\in
\Theta$ by formula (\ref{liml0})
is at least three times continuously differentiable on compact $\overline{\Theta}$ by {\sc (H3b)}, and
\[
\begin{array}{l}
D\ell_{\vartheta_0} (\vartheta )=\int_E \frac{
(\mu (x,\vartheta_0 )-\mu (x,\vartheta ))}{b^2(x)}D_\vartheta \mu (x,\vartheta )\,\pi_{\vartheta_0} (dx)\\
D^2\ell_{\vartheta_0} (\vartheta )=\int_E\left(\frac{
(\mu (x,\vartheta_0 )-\mu (x,\vartheta ))}{b^2(x)}D_\vartheta^2 \mu (x,\vartheta )
-\frac{1}{b^2 (x)}(D_\vartheta^\tau \mu D_\vartheta\mu )(x,\vartheta)\right)\,\pi_{\vartheta_0} (dx).
\end{array}
\]
Hence, by the same argument as for (\ref{liml0}),
for any fixed $\vartheta\in\overline{\Theta}$, $\Pb_{(\vartheta_0,\sigma)}$-a.s.
\begin{equation}\begin{array}{lcl}
\lim_{T\rightarrow +\infty}\frac{1}{T}D\ell_T (\vartheta ) &=& D\ell_{\vartheta_0} (\vartheta ),\\
\lim_{T\rightarrow +\infty}\frac{1}{T}D^2\ell_T (\vartheta ) &=& D^2\ell_{\vartheta_0} (\vartheta ).
\end{array}\label{liml0dd2}
\end{equation}
If $\vartheta\neq\vartheta_0$ then $\ell_{\vartheta_0}(\vartheta )<\ell_{\vartheta_0} (\vartheta_0 )$
by (\ref{liml0}), 
and 
{\sc (H4b)}. Hence $\vartheta_0$ is the unique point of maximum of
$\ell_{\vartheta_0}$ on $\overline{\Theta}$. This implies identifiability property of the model: let $\vartheta_1,\vartheta_2\in\Theta$ be such that $\Pb_{(\vartheta_1,\sigma )}=\Pb_{(\vartheta_2,\sigma )}$. Then
$\pi_{\vartheta_1}=\pi_{\vartheta_2}$ and so $\ell_{\vartheta_1}\equiv\ell_{\vartheta_2}$ by (\ref{liml0}). Hence
$\vartheta_1 =\vartheta_2$. Moreover, {\sc (H5b)} implies that the Fisher information matrix is positive definite,  i.e.\
\[\begin{array}{l}
 I(\vartheta_0 )=-D^2\ell_{\vartheta_0} (\vartheta_0 )=\int_E\frac{1}{b^2 (x)}(D_\vartheta^\tau \mu D_\vartheta\mu )(x,\vartheta_0 )\,\pi_{\vartheta_0} (dx)>\O.
\end{array}\]
The next theorem states that the continuous-time MLE of drift parameters exists, is consistent and asymptotically efficient, and satisfies assumptions {\sc (H4a)} and {\sc (H6a)} a.s.\ for almost all observational times. Generally
these are well known facts (see e.g. \cite{BH} or \cite{Feigin}) but we provided it here for completeness, and in the appropriate form for the purpose of proving Theorem \ref{tm:erg} below.

\begin{theorem}\label{tm:cont}
Let us assume that {\sc (H1b-5b)} hold.  Then
there exists an $({\cal F}_T^0, T>0)$-adapted process $(\hat{\vartheta}_{T}, T>0 )$ of random vectors
such that for every $\theta =(\vartheta ,\sigma )\in\Psi$ the following holds:
\begin{romannum}
\item $\Pb_{\theta}$-a.s.\ there exists $T_0 >0$ such that
for all $T\geq T_0$, $\hat{\vartheta}_{T}\in\Theta$ is the unique point of maximum of  $\ell_T$ on $\overline{\Theta}$, and $D^2 \ell_T (\hat{\vartheta}_T)<\O$ in a way that $\min_{|y|=1}y^\tau (-\frac{1}{T}D^2\ell_T(\hat{\vartheta}_T))y\geq\frac{1}{2}\min_{|y|=1}y^\tau I(\vartheta )y$.
\item  $\lim_{T\rightarrow +\infty}\hat{\vartheta}_T=\vartheta$ $\Pb_\theta$-a.s.
\item  $(\sqrt{T}(\hat{\vartheta}_T-\vartheta ), T>0)$ converges in law w.r.t.\ $\Pb_{\theta}$ to normal law $N(\O,\sigma I(\vartheta)^{-1})$ with expectation $\O$  and covariance matrix
    $\sigma I(\vartheta)^{-1}$.
\end{romannum}
\end{theorem}\vspace{5pt}

\noindent
The following theorem is a version of Theorem \ref{tm:30} for ergodic diffusions. In addition it states that AMLEs are consistent and asymptotically efficient when both maximal observational time and number of discrete observational time points tend to infinity for appropriate sampling schemes. Hence in its statement '$\lim_{n,T}$' denotes the limit when
both $T\rightarrow +\infty$ and $n\rightarrow +\infty$.

\begin{theorem}\label{tm:erg}
Let us assume that {\sc (H1b-5b)} hold.
Then there exists a process $(\hat{\vartheta}_{n,T}; n\geq 1, T>0 )$ of ${\cal F}_{n,T}$-measurable random vectors
$\hat{\vartheta}_{n,T}$ such that for all $\theta =(\vartheta ,\sigma )\in\Psi$ and $\pi_\vartheta$-a.s.\ nonrandom initial conditions,
and all equidistant samplings such that $\delta_{n,T}=T/n\rightarrow 0$, the following holds.
\begin{romannum}
\item $\lim_{n,T}\Pb_{\theta}(D{\ell}_{n,T}
(\hat{\vartheta}_{n,T})=\O )=1$.
\item $(\Pb_{\theta})\lim_{n,T}(\hat{\vartheta}_{n,T}-
\hat{\vartheta}_T)=\O$,
\item $\hat{\vartheta}_{n,T} -\hat{\vartheta}_T=O_{\Pb_\theta}(\sqrt{\delta_{n,T}})$, $n\rightarrow +\infty$, $T\rightarrow +\infty$
 \item If $(\tilde{\vartheta}_{n,T};n\geq 1,T>0)$ is a process of random vectors in $\Theta$ that satisfies $(i-ii)$ then
 $\lim_{n,T}\Pb_{\theta} (\tilde{\vartheta}_{n,T} =\hat{\vartheta}_{n,T})=1$.
 \item $(\Pb_{\theta})\lim_{n,T}\hat{\vartheta}_{n,T}=
\vartheta$, and if
in addition $\lim_{n,T}T\delta_{n,T}=0$ then
\[\begin{array}{l}
\sqrt{T}(\hat{\vartheta}_{T,n}-\vartheta )\stackrel{{\cal L}-\Pb_{\theta}}{\longrightarrow} N(\O,\sigma I(\vartheta)^{-1}),\; T\rightarrow +\infty, n\rightarrow +\infty.
\end{array}\]
\item $(\Pb_{\theta})\lim_{n,T}\sigma_{n,T}=\sigma$, and if in addition $\lim_{n,T}T\delta_{n,T}=0$ then
\[\begin{array}{l}\sqrt{n}\frac{1}{\sigma\sqrt{2}
}(\hat{\sigma}_{n,T} -\sigma)\stackrel{{\cal L}-\Pb_{\theta}}{\longrightarrow} N(0,1),\;
T\rightarrow +\infty,n\rightarrow +\infty.\end{array}\]
\end{romannum}
\end{theorem}

\section{Examples}

\begin{example}\label{primjer1} Generalized logistic model. {\rm  Let the stochastic generalized logistic model be given with the following SDE:
\begin{equation}\begin{array}{l}
dX_t =(\alpha-\beta X_t^\gamma)X_t\, dt + \sqrt{\sigma} X_t\, dW_t,\; X_0 =x_0>0\label{logSDE}
\end{array}\end{equation}
where $\vartheta =(\alpha,\beta,\gamma)$ ($\gamma>0$) is a drift vector parameter. By using the methods of
stochastic calculus it is possible to explicitly solve (\ref{logSDE}) that proves
that there exists pathwise unique, continuous and strong solution to this SDE
with $X$ defined on $\Omega\times [0,+\infty\rangle$ and values in $E=\langle 0,+\infty\rangle$.
Moreover, it turns out that for drift parameters such that
$\alpha>\sigma/2$, $\beta>0$ and $\gamma >0$, generalized logistic process $X$ is positive recurrent and ergodic
with a such stationary distribution $\pi_\vartheta$ that for stationary $X$, $X_t^\gamma$ follows
$\Gamma$-distribution with parameters $A:=2(\alpha-\sigma/2)/(\gamma\sigma)$ and $B:=\gamma\sigma/(2\beta)$ (i.e.
$\E X_t^\gamma =AB$, $\E (X_t^{\gamma})^2 =AB(B+1)$) by e.g. Theorem 7.1, pp. 219-220 in \cite{Friedman}.
Hence, assumption {\sc (H1b)} holds.

In generalized logistic model, drift function is equal $\mu (x,\vartheta )=(\alpha-\beta x^\gamma )x$, and up to the
diffusion parameter $\sigma >0$, diffusion coefficient function is $b(x)=x >0$ on $E$. Hence $b'\equiv 1$,
$bb''=b^2 b'''\equiv 0$ that are trivially integrable with respect to any probability law. Let
$f (x,\vartheta )=\mu (x,\vartheta)/b(x)=\alpha-\beta x^\gamma $. Notice that any partial derivatives of $f$ with respect
to $\vartheta$ are of the form $-\beta^n x^\gamma \log^m x$ where $n\in\{ 0,1\}$, $m\in\N_0$. Of the same forms
are components of $b^{k} \frac{\partial^k}{\partial x^k}D_\vartheta^m f$ for $k=1, 2$. Finally, any $p$-th power of
their absolute values ($p$ is a positive integer) are of the form $x^c |\log x|^m$ up to a constant, where $c >0$ is a real
number and $m$ is a nonnegative integer. These functions are
integrable with respect to $\pi_\vartheta$. If we choose a relative compact $\Theta$ of drift parametric set
$\langle\frac{\sigma}{2},+\infty\rangle\times\langle 0,+\infty\rangle^2$ then there exist $\alpha_0>\sigma/2$,
$\beta_0 >0$ and $\gamma_0>0$ such that for all $\vartheta\in\overline{\Theta}$, $x>0$, and all $0\leq m\leq 6$, $k\in\{0,1,2\}$
and integers $j_\alpha$, $j_\beta$, $j_\gamma$ such that $j_\alpha+j_\beta +j_\gamma=m$,
\[\begin{array}{l}
|b^k(x) \frac{\partial^{m+k}}{\partial^k x\partial\alpha^{j_\alpha}\partial\beta^{j_\beta}\partial\gamma^{j_\gamma}} f(x,\vartheta)|\leq g(x):=\alpha_0+\beta_0 x^{\gamma_0}(1+\log^2 x+\log^4 x +\log^6 x).
\end{array}\]
Then $g\in L^p (\pi_\vartheta)\cap C^1 (E)$ for all $p\geq 1$ and $\vartheta\in\overline{\Theta}$ which implies
partially {\sc (H2b)} and {\sc (H3b)}  by simple calculation (see the proof of Corollary \ref{tm:cor2} below). To finish the proof of {\sc (H2b)} notice that for all $h_0>0$, and all $0<h\leq h_0$,
\[\begin{array}{l}
\exp(16\int_0^h((\alpha -\beta X_t^\gamma+\frac{15\sigma}{2})\, dt)\leq \exp((16\alpha_0+120\sigma) h_0)=c(x_0)=\mbox{\rm constant}
\end{array}\]
since $X_t >0$ for all $t\geq 0$ and $\beta >0$. This implies the same inequalities for expectations with respect to any
initial conditions $X_0 =x_0$. Hence {\sc (H2b)} is proved.

To show that {\sc (H4b)} holds, let us assume that
\[\begin{array}{l}\int_E (\mu (x,\vartheta_1)-\mu (x,\vartheta_2))^2/b^2(x)\pi_{\vartheta_1}(dx)=0
\end{array}\]
for some $\vartheta_1\in\Theta$ and $\vartheta_2\in\overline{\Theta}$. Since $\pi_{\vartheta_1}$ is absolutely continuous w.r.t.
Lesbegues measure $\lambda$ on $E$, this implies that $\mu(x,\vartheta_1)=\mu (x,\vartheta_2)$ for a.s. $x>0$ w.r.t. $\lambda$.
Hence, smooth function $u(x):=\beta_1 x^{\gamma_1}-\beta_2 x^{\gamma_2}$ must be a constant function for $\lambda$-a.s. $x>0$. This
implies that $\gamma_1=\gamma_2$ and hence $\vartheta_1 =\vartheta_2$. This proves  {\sc (H4b)}.

Finally, {\sc (H5b)} holds since $\frac{\partial}{\partial\alpha} \mu (x, \vartheta)=1$,  $\frac{\partial}{\partial\beta} \mu (x, \vartheta)=-x^\gamma$, and
$\frac{\partial}{\partial\gamma} \mu (x, \vartheta)=-\beta x^\gamma\log x$ are obviously linearly independent functions in $L^2(\pi_\vartheta )$.}
\end{example}

\begin{example}\label{primjer2} Cox-Ingersoll-Ross (CIR) model. {\rm CIR model (or Feller's square root model) is given by SDE:
\begin{equation}\label{cirSDE}
\begin{array}{l}
dX_t =(\beta-\alpha X_t)\, dt + \sqrt{\sigma |X_t |}\, dW_t,\; X_0 =x_0>0.
\end{array}
\end{equation}
Vector of drift parameters is $\vartheta =(\alpha,\beta)$, drift function $\mu (x,\vartheta )=\beta -\alpha x$ is linear in its parameters, and $b(x)=\sqrt{|x|}$. It has been known (see e.g. \cite{Jin}) that if $\alpha >0$ and $\beta >0$ are such that $2\beta >\sigma$, and $x_0>0$ then SDE (\ref{cirSDE}) has strong positive recurrent and ergodic solution in state space $E=\langle 0,+\infty\rangle$ with stationary distribution $\pi_\vartheta$ which has $\Gamma$-law  with
expectation $\beta/\alpha$ and variance $\beta \sigma/(2\alpha^2)$. Hence {\sc (H1b)} and {\sc (H2a-3a)} hold for any open relatively compact and convex set $\Theta$ in $\langle 0,+\infty\rangle^2\cap\{ (\alpha,\beta ):2\beta >\sigma\}$ that contains the true drift parameter value. Additionally let us assume that if $\vartheta =(\alpha,\beta )\in\overline{\Theta}$ then $2\beta/\sigma >16$. Then function $x\mapsto 1/x$ is in $L^{16}(\pi_\vartheta )$ which implies {\sc (H3b)} and partially {\sc (H2b)}. Since inequality in {\sc (H2b)} is used only for proving the statement of Lemma \ref{lema1} it is sufficient to prove this lemma directly (instead of this inequality), i.e. for each
$\vartheta\in\Theta$ we want to find
a function $c_0\in L^1 (\pi_\vartheta )$ and $h_0>0$ such that the following inequality holds for any $t\geq 0$:
\begin{equation}\label{cirineq}\begin{array}{l}
\sup_{0<h\leq h_0}\E_{(\vartheta,\sigma)}(b(X_{t+h}/b(X_t ))^8\leq\E_{(\vartheta,\sigma)}c_0(X_t).
\end{array}
\end{equation}
Let $x >0$ and $h>0$ be arbitrary, and let $X$ be such that (\ref{cirSDE}) holds with $X_0 =x$. Let $\E\equiv\E_{(\vartheta,\sigma)}$,
and let us calculate
\[\begin{array}{l}
\E (b(X_{h})/b(x ))^8=(1/x)^4\cdot\int_0^{+\infty} y^4 p(h,x,y)\, dy
\end{array}
\]
where
\[\begin{array}{l}
p(h,x,y)=C e^{-u-v}(v/u)^{q/2} I_q(2\sqrt{uv}),
\end{array}
\]
is the transition density of $X_h$ given $X_0 =x$ (see \cite{Jin}). Here $u=Cxe^{-\alpha h}$, $v=Cy$, $C=(2\alpha)/(\sigma (1-e^{-\alpha h}))$, $q=(2\beta/\sigma)-1$, and $I_q$ is the modified Bessel function of the first kind of order $q$. Since
$q >15$ and $e^{-\alpha h}<1$ it turns out that
\[\begin{array}{l}
\E (b(X_{h})/b(x ))^8=(1/x)^4\cdot\int_0^{+\infty} y^4 p(h,x,y)\, dy\leq 8q^4 (\frac{3}{x^4}+\frac{12}{x^3}+
\frac{9}{x^2}+\frac{2}{x}+1)=:c_0(x).
\end{array}
\]
Then
\[\begin{array}{l}
\E (b(X_{t+h}/b(X_t ))^8=\E [\E [(b(X_{t+h}/b(X_t ))^8|{\cal F}_t^0]]=\E [\E_{X_t}[(b(X_{h}/b(X_0 ))^8]]\leq\E c_0(X_t)
\end{array}
\]
by Markov property and above inequality. Hence (\ref{cirineq}) holds for any $h_0>0$, and $c_0\in L^1 (\pi_\vartheta )$ since $x\mapsto 1/x\in L^{16}(\pi_\vartheta )$.

Finally, {\sc (H4b-5b)} follow easily since functions $\sqrt{x}$ and $1/\sqrt{x}$ are linearly independent, and $\pi_\vartheta$ is dominated by Lesbegue measure on $\langle 0,+\infty\rangle$. Hence, if $2\beta/\sigma >16$ then
Theorem \ref{tm:erg} can be applied on CIR model (\ref{cirSDE}).

Since in CIR model the drift function is linear in its parameters ALF $\ell_{n,T} (\vartheta )$  and LF $\ell_T (\vartheta )$ ($\vartheta =(\alpha,\beta )$) are quadratic functions. Hence there exist unique explicit solutions to stationary equations $D\ell_{n,T}(\vartheta )=0$ and $D\ell_{T}(\vartheta )=0$, and properties of the AMLE can be investigate by simulation techniques easily. For this purpose we simulate $M=1000$ paths of the process $X$ over
time-interval $[0,T]$ for true parameter values $\vartheta_0 =(\alpha_0,\beta_0)=(0.5,0.03)$ and $\sigma_0 =0.06^2$,
and several different values of $T$, precisely for $T=3,4,\ldots,11$. Drift parameter values $\vartheta_0$ have
been borrowed from similar examples in \cite{AitSahalia0} or \cite{Li}, and $\sigma_0$ has been
chosen to be a such that $2\beta_0/\sigma_0 \approx 16.7>16$. Each path initially starts at $x_0=1$, and have been simulated by using Milstein sheme based on discretization of $[0,T]$ on $2^{16}$ equidistant points. Using the same discretization $[0,T]$ each Riemann integral in $\ell_T (\vartheta )$ have been approximated by trapezoidal rule, and It\^{o} integral by Euler approximation. Any
estimate $\hat{\theta}_{n,T}=(\hat{\vartheta}_{n,T},\hat{\sigma}_{n,T})$ for varying $n$ has been calculated from the same path as estimate $\hat{\vartheta}_T$ does.

The results of analyzing asymptotic behavior of deviances $\hat{\vartheta}_{n,T}-\hat{\vartheta}_T$ and
$\hat{\sigma}_{n,T}-\sigma_0$ are presented at Figure 1. 
\begin{figure}[p]
\hfil \scalebox{0.9}{\includegraphics{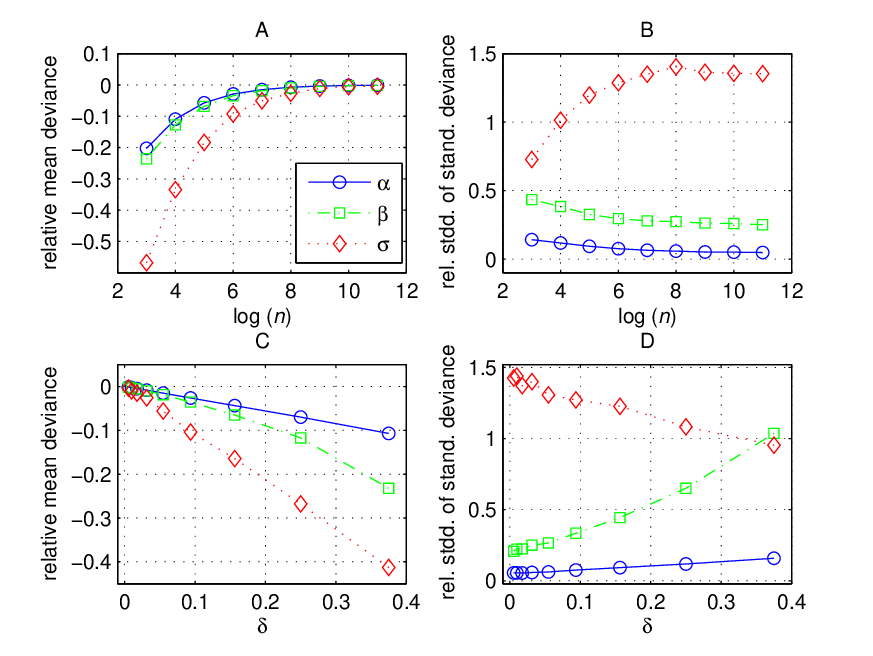}}\hfil
\caption{({\sc A}) Relative means of components of statistics $\hat{\vartheta}_{n,T}-\hat{\vartheta}_T$ and
$\hat{\sigma}_{n,T}-\sigma_0$ (relative to $\vartheta_0$ and $\sigma_0$ respectively) with respect to different
sampling sizes $n$ and fixed $T=7$. ({\sc B}) Relative standard deviations of standardized deviances $\sqrt{T}(\hat{\vartheta}_{n,T}-\hat{\vartheta}_T)/\sqrt{\delta_{n,T}}$ and $\sqrt{n}(\hat{\sigma}_{n,T}-\sigma_0)/\sigma_0$ (relative to the true parameter values) with respect to different
sampling sizes $n$ and fixed $T=7$. ({\sc C}) Relative means of components of the same statistics as in {\sc A} but
with respect to $\delta=\delta_{n,T}=T/2^T =\log (n)/n$
for different $T$s. ({\sc D}) Relative standard deviations of the same standardized deviances as in {\sc B} but
with respect to $\delta=\delta_{n,T}=T/2^T =\log (n)/n$ for different $T$s. In all cases means and std. deviatians are estimated based on simulated samples with length $M=1000$.}
\end{figure}
Subfigures {\sc A} and {\sc C} represent mean deviances relative to the true parameter values, and subfigures {\sc B} and {\sc D} represents standard deviations of deviances standardized with $\sqrt{\delta_{n,T}}/\sqrt{T}=1/\sqrt{n}$ and relative to the true parameter values too.

In case of subfigures {\sc A} and {\sc B}, $T=7$ is fixed, but
number $n$ of equidistant sampling time-points varies from $2^3$ to $2^{11}$ in a way that $\log (n)=k$, $k=3,4,\ldots,11$, where '$\log (\cdot )$' represents logarithm with base 2.
Subfigure {\sc A} shows the expected asymptotic behavior that $\lim_n\hat{\vartheta}_{n,T}=\hat{\vartheta}_T$
and $\lim_n\hat{\sigma}_{n,T}=\sigma_0$ in case of fixed $T$ and $\delta_{n,T}=T/n\rightarrow 0$, but also that
AMLE subestimates MLE and similarly for $\hat{\sigma}_{n,,T}$.
The rate of convergence can be seen from subfigure {\sc B}. Namely, the convergence of empirical standard deviations (estd)
of components of $\sqrt{T}(\hat{\vartheta}_{n,T}-\hat{\vartheta}_T)/\sqrt{\delta_{n,T}}$ (relative to $\vartheta_0$) shows that these statistics are bounded in probability, while the convergence of estd of $\sqrt{n}(\hat{\sigma}_{n,T}-\sigma_0)/\sigma_0$ to a neighborhood of $\sqrt{2}\approx 1.41$ are also expected by convergence in law of $\sqrt{n}(\hat{\sigma}_{n,T}-\sigma_0)$ to the normal distribution with standard deviation $\sigma_0\sqrt{2}$. Table 1 
shows p-values of three tests of normality: Shapiro-Wilk (SW), Lilliefors Kolmogorov-Smirnov (Lillie) and Jarque-Bera (JB), and Kolmogorov-Smirnov test (KS) of standard normality of simulated statistic $(\sqrt{n}/\sigma_0\sqrt{2})(\hat{\sigma}_{n,T}-\sigma_0)$ with respect to $n$ (and fixed $T=7$).
\begin{table}
\[\begin{array}{crrrr}\hline
 \log (n) & \mbox{\rm SW} & \mbox{\rm Lillie} & \mbox{\rm JB} & \mbox{\rm KS}\\ \hline
 3 & 0.0000  &   0.0010  &  0.0010  &  0.0000\\
 4 & 0.0000  &  <0.0010  & <0.0010  &  0.0000\\
 5 & 0.0000  &  <0.0010  & <0.0010  &  0.0000\\
 6 & 0.0001  &   0.0136  & <0.0010  &  0.0000\\
 7 & 0.0190  &   0.0382  &  0.0240  &  0.0000\\
 8 & 0.0959  &   0.1104  &  0.0926  &  0.0000\\
 9 & 0.3887  &  >0.5000  &  0.2146  &  0.0000\\
 10 &   0.5968 &   >0.5000  &  >0.5000 &   0.0010\\
 11 &  0.6537  &  >0.5000   & >0.5000   & 0.0037\\ \hline
 \end{array}\]
 \caption{P-values of Shapiro-Wilk (SW), Lilliefors (Lillie), Jarque-Bera (JB) and Kolmogorov-Smirnov (KS) tests
 of normality applied on samples of
 statistic $(\sqrt{n}/\sigma_0\sqrt{2})(\hat{\sigma}_{n,T}-\sigma_0 )$ (of length $M=1000$) with respect to different sampling sizes $n$ with fixed $T=7$.}
\end{table}
Obviously, the statistic converges to normality, but slowly to the specific limiting normal distribution.

The same behavior of deviances $\hat{\vartheta}_{n,T}-\hat{\vartheta}_T$ and
$\hat{\sigma}_{n,T}-\sigma_0$ when $T\rightarrow +\infty$ in a way that $\delta_{n,T}\rightarrow 0$ can be seen
from subfigures {\sc C} and {\sc D}. In case of these subfigures the relative mean deviances and the relative
standard deviations of standardized deviances are presented with respect to $\delta=\delta_{n,T}=T/2^T =\log (n)/n$
for $T=\log (n)=3,4\ldots,11$. Normal q-q plot of the sample of $(\sqrt{n}/\sigma_0\sqrt{2})(\hat{\sigma}_{n,T}-\sigma_0)$ in case $T=11$ and $n=2^{11}$ is
presented at subfigure B of Figure 2.

Asymptotic properties of deviances $\hat{\vartheta}_{n,T}-\vartheta_0$ when $T\rightarrow +\infty$ in a way that $\delta_{n,T}\rightarrow 0$, are presented in Figure 2. Subfigure {\sc A} presents the relative mean deviances with
respect to $\delta=T/2^T =\log (n)/n$ for $T=\log (n)=3,4\ldots,11$. We notice the concave shape of the both curves tending to zero when $\delta\rightarrow 0$. The convergence to the normality is very slow as illustrated with
q-q plots of the standardized components of AMLEs (with respect to the limiting normal laws and at $T=11$, $n=2^{11}$) at subfigures {\sc C} and {\sc D}.
\begin{figure}[p]
\hfil \scalebox{0.9}{\includegraphics{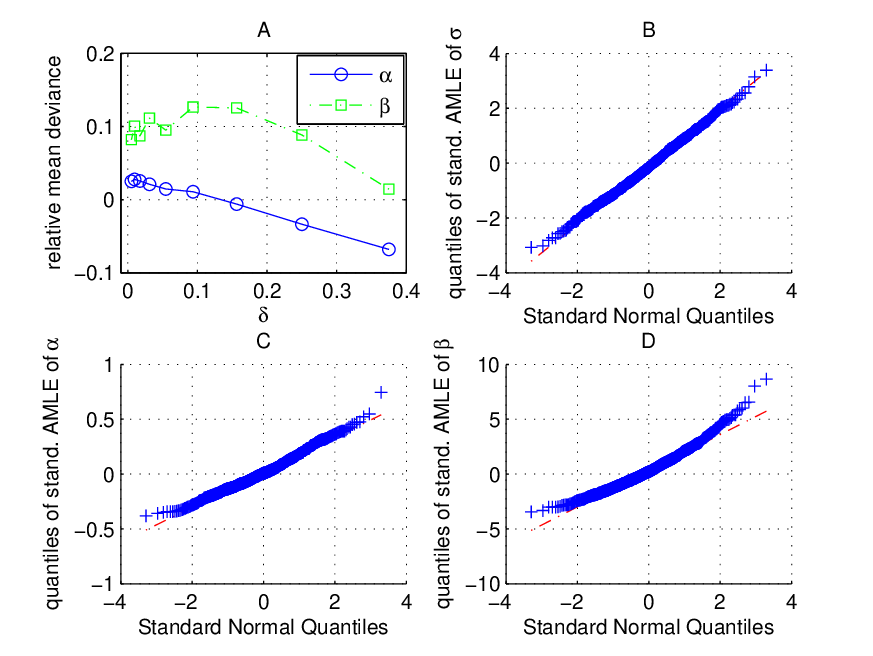}}\hfil
\caption{({\sc A}) Relative means of components of statistics $\hat{\vartheta}_{n,T}-\vartheta_0$ (relative to $\vartheta_0$) with respect to $\delta=\delta_{n,T}=T/2^T =\log (n)/n$
for different $T$s. ({\sc B}) Normal q-q plot of the standardized statistics $\sqrt{n}(\hat{\sigma}_{n,T}-\sigma_0)$ with respect to the limiting normal law (with $T=11$, $n=2^{11}$) ({\sc C-D}) Normal q-q plots of the standardized components of statistics $\sqrt{T}(\hat{\vartheta}_{n,T}-\vartheta_0)$ with respect to the limiting normal law (with $T=11$, $n=2^{11}$) ({\sc C} for $\alpha$ and {\sc D} for $\beta$ components).
In all cases estimations are based on simulated samples with length $M=1000$.
}
\end{figure}
}
\end{example}

\section{Proofs}

Basically the
proof of Theorem \ref{tm:30} is based on the so-called general
theorem on approximate maximum likelihood estimation and its
corollary that are stated and proved in \cite{huzakMLE} as Theorem 3.1
and Corollary 3.2. The proof of Theorem \ref{tm:erg} is a modification of
the proof of the same theorem based on Theorem \ref{tm:cont}.
But first we need to state and prove Theorem \ref{tm:new}, and
its Corollaries \ref{tm:23} and \ref{tm:cor2} that are needed in applying the general theorem in
this context. Proofs of some technical lemmas are in Appendix.

Let us suppose that $X=(X_t, t\geq 0)$ is a diffusion satisfying {\sc (H1a-2a)} with true parameter
$\theta_0 =(\vartheta_0,\sigma )\in\Theta$, and such that $\Pb (X_0=x_0)=1$ for $x_0\in E$. Here
$\Pb\equiv\Pb_{\theta_0}$, $\E\equiv\E_{\theta_0}$ and $L^2\equiv L^2(\Pb )$.  Denote $\mu_0 =\mu (\cdot,\vartheta_0)$, $\nu =\sqrt{\sigma} b$,
${\cal A}a:=a'\mu_0+a''\nu^2/2$, and $\bar{a}:=|{\cal A}a|+|a'\nu|$ for $a\in C^2(E)$.

\begin{theorem}\label{tm:new}
Let $\Theta\subset\R^d$ be an open convex set, and let $f:E\times\Theta\rightarrow\R$,
$a:E\rightarrow\R$ be functions. Let $0=t_0<t_1<\cdots <t_n=T$ be subdivisions of intervals $[0,T]$, $T>0$, such that
$\delta_{n,T}\downarrow 0$. Assume the following:
\begin{nitemizetm}
\item[{\sc (B1):}] $a\in C^2(E)$ and there exist constants $C_a>0$, $T_a\geq 0$, and $n_a\in\N$ such that
\[\begin{array}{l}
(\forall\, T >T_a)(\forall n\geq n_a)\;\;\frac{1}{T}\E\left(\!\!\int_{0}^{T}\!\!(a^4\!+\!\bar{a}^4)(X_t))dt
+\sum_{i=0}^{n-1}\!a^4 (X_{t_i})\Delta_i t\right)\leq C_a.
\end{array}\]
\item[{\sc (B2):}] For all $\vartheta\in\Theta$, $f(\cdot,\vartheta )\in C^2 (E)$, and for all
$(x,\vartheta)\in E\times\Theta$ and $1\leq m\leq d+1$ there exists partial derivatives
$D^m_\vartheta f(x,\vartheta )$, $\frac{\partial}{\partial x}D^m_\vartheta f(x,\vartheta )$, and
$\frac{\partial^2}{\partial x^2}D^m_\vartheta f(x,\vartheta )$.  Moreover,
\[\begin{array}{l}
(\forall\, 0\leq m\leq d+1)\;\;
D^m_\vartheta f,\, \frac{\partial}{\partial x}D^m_\vartheta f,\, \frac{\partial^2}{\partial x^2}D^m_\vartheta f\in C(E\times\Theta).
\end{array}\]
\item[{\sc (B3):}] For any relatively compact set ${\cal K}$ in $\Theta$ there exist: a positive measurable function
$g:E\rightarrow\R$ such that for all $0\leq m\leq d+1$,
\[\begin{array}{l}
\sup_{\vartheta\in\overline{\cal K}}\left(\!|D^m_\vartheta f (\cdot ,\vartheta )|_\infty\!+\!|\frac{\partial}{\partial x}D^m_\vartheta f (\cdot ,\vartheta )|_\infty(|\mu_0 |\! +\! |\nu |)\!+\!
|\frac{\partial^2}{\partial x^2}D^m_\vartheta f (\cdot ,\vartheta  )\nu^2 |_\infty\right)\! \leq\! g,
\end{array}\]
and constants $C_g>0$, $T_g\geq 0$, and $n_g\in\N$, such that
\[\begin{array}{ll}
(\forall\, T>T_g)(\forall n\geq n_g)\;\;\frac{1}{T}\E\left(\!
\int_{0}^{T}g^4\! (X_t)\,dt +
\sum_{i=0}^{n-1}g^4\! (X_{t_i})\Delta_i t\!\right) & \leq C_g\;\&\\
\quad\quad\frac{1}{T}\E 
\sum_{i=0}^{n-1}\left(\! (ga)^4\! (X_{t_i})\Delta_i t +
g^4\!(X_{t_i})\int_{t_i}^{t_{i+1}}(a^4+\bar{a}^4)\!(X_t)\,dt\!\right) &\leq C_{g}.
\end{array}\]
\item[{\sc (B4):}] There exist: a measurable function $c: E\rightarrow\R$ and constants $h_0>0$, $C_c>0$, $T_c\geq 0$, and $n_c\in\N$
 such that for $r:=|\frac{\mu_0 b'}{b}| +|b''b|+|b'|$, 
\[\begin{array}{ll}
\sup_{0<h\leq h_0}\E\exp\left( 8\int_0^{h}\left(2\frac{\mu_0 b'}{b} +\sigma (b''b+15b'^2)\right)(X_s)\,ds\right) &\leq c(x_0),\\
(\forall\, T>T_c)(\forall n\geq n_c)\;\;\frac{1}{T}\E\left(\!
\sum_{i=0}^{n-1}c(X_{t_i})\Delta_i t + 
\int_0^T r^8(X_t)dt\right) &\leq C_c.
\end{array}\]
\end{nitemizetm}
Then there exist constants $C_1>0$, $C_2>0$, $T_0 \geq 0$, and $n_0\in\N$, possible dependent on ${\cal K}$, $d$, and $a$, such that
for all $T> T_0$, and $n\geq n_0$, 
\begin{eqnarray}
\E\sup_{\vartheta\in\overline{\cal K}}\left(\frac{1}{T\sqrt{\delta_{n,T}}}\sum_{i=0}^{n-1}\int_{t_i}^{t_{i+1}}\!\!\!(f(X_t,\vartheta ) -f(X_{t_i},\vartheta))a(X_t)\,dt\right)^2\!
&\leq & \! C_1\label{ineq1}\\
\E\sup_{\vartheta\in\overline{\cal K}}\left(\frac{1}{\sqrt{T\delta_{n,T}}}\sum_{i=0}^{n-1}\int_{t_i}^{t_{i+1}}\!\!\!(f(X_t,\vartheta ) -f(X_{t_i},\vartheta))a(X_t)\,dW_t \right)^2\!
&\leq & \! C_1\label{ineq2}\\
\E\sup_{\vartheta\in\overline{\cal K}}\left(\frac{1}{T\sqrt{\delta_{n,T}}}\sum_{i=0}^{n-1}\int_{t_i}^{t_{i+1}}\!\!\!f(X_{t_i},\vartheta )\left(\frac{b(X_t )}{b(X_{t_i})}-1\right) a(X_t)\,dt \right)^2\! &\leq &\! C_2 \label{b2}\\
\E\sup_{\vartheta\in\overline{\cal K}}\left(\frac{1}{\sqrt{T\delta_{n,T}}}\sum_{i=0}^{n-1}\int_{t_i}^{t_{i+1}}\!\!\!f(X_{t_i},\vartheta )\left(\frac{b(X_t )}{b(X_{t_i})}-1\right) a(X_t)\,dW_t \right)^2\! &\leq &\! C_2. \label{b3}
\end{eqnarray}
\end{theorem}

\begin{remark}\label{rem3}{\rm If function $f$ that satisfies (B2) and all its partial derivatives from (B2) are bounded on $E\times\overline{\cal K}$ then (B3) holds if $a$ is bounded too.
Similarly if $a\in C^2(E)$
is bounded then satisfies (B1).
 If in addition $\mu_0$, $b$, $b'$, and $b''$ are bounded then (B4) holds for constant function $c\equiv\exp(\gamma h_0)$ where $\gamma >0$ and $h_0>0$ are constants. In this case the statements of Theorem \ref{tm:new} hold for $T_0=0$, and hence for
all $T>0$  obviously from the proof of Theorem \ref{tm:new}.
}
\end{remark}

For a moment let us assume that ${\cal K}=\prod_{i=1}^d\langle a_i,b_i\rangle$ is an open and bounded
$d$-dimensional rectangular in $\Theta$. Then there
exists $\varepsilon >0$ such that ${\cal K}_\varepsilon:=\prod_{i=1}^d\langle a_i-\varepsilon,b_i+\varepsilon\rangle$ is an open and bounded $d$-rectangular in $\Theta$ too. Let $\phi:\R^d\rightarrow\R$ be a $C^\infty$-function such that
$\phi\equiv 1$ on ${\cal K}$ and $\phi\equiv 0$ on ${\cal K}_\varepsilon^c$. Such a function exists (see e.g.\ \cite{Borisovich}, Lemma IV.4.4, p.\ 176). Then function $(x,\vartheta )\mapsto \tilde{f}(x,\vartheta):=f(x,\vartheta )\cdot\phi (\vartheta )$
satisfies (B2-3) if $f$ satisfies the same assumption (with rescaled function $g$). Namely, $\tilde{f}\equiv f$ on $E\times{\cal K}$
and $\tilde{f}\equiv 0$ on $\partial{\cal K}_\varepsilon$. The same holds for all partial derivatives of $\tilde{f}$ that exist, and
$\tilde{f}$ satisfies (B1) obviously. Since
$\phi$ and all of its derivatives are bounded, $\tilde{f}$ satisfies (B3) too with $Cg$ instead of $g$ with a constant $C$ depending on $\phi$.
Obviously, statements (\ref{ineq1}-\ref{b3}) hold for a function $f$ that satisfies (B2-3), and a rectangular
${\cal K}$
if (\ref{ineq1}-\ref{b3}) hold for $\tilde{f}$ and the rectangular ${\cal K}_\varepsilon$. Moreover, notice that if (\ref{ineq1}-\ref{b3}) hold for an arbitrary open and bounded $d$-dimensional rectangular ${\cal K}$, then the same statements
hold for every relatively compact set in $\Theta$. Hence it is sufficient to prove (\ref{ineq1}-\ref{b3}) for an open and bounded $d$-dimensional rectangular ${\cal K}\subset\Theta$, and a function $f$ satisfying (B2-3) and the following additional assumption.\vspace{5pt}

{\sc (B ${\cal K}$)}: For all $x\in E$ and all $0\leq m\leq d+1$,
$D^m_\vartheta f(x,\cdot)\equiv\O$, $\frac{\partial}{\partial x}D^m_\vartheta f(x,\cdot )\equiv\O$
and $\frac{\partial^2}{\partial x^2}D^m_\vartheta f(x,\cdot )\equiv\O$
on $\partial{\cal K}$.\vspace{5pt}

Moreover, let $A$ be an invertible affine mapping of $\R^d$, and let $f$ be a function
on $E\times\Theta$ that satisfies (B2-3) and (B ${\cal K}$). Then the function $\bar{f}$ defined on
$E\times A(\Theta )$ by the rule $\bar{f}(x,\eta ):=
f(x,A^{-1}\eta )$, satisfies (B2-3) and (B $A{\cal K}$) too. Since the left hand side of (\ref{ineq1}-\ref{b3})
do not change by the change of variable $\vartheta\mapsto \eta=A\vartheta$, it is sufficient
to prove (\ref{ineq1}-\ref{b3}) for ${\cal K}_0:=\langle -\pi,\pi\rangle^d$ and a function $f$ that satisfies
(B2-3) and (B ${\cal K}_0$).

Now, let $f$ be a function satisfying (B2-3) and (B ${\cal K}_0$). For $x\in E$,
$\vk=(k_1,\ldots,k_d)\in\Z^d$, $\vartheta =(\vartheta_1,\ldots,\vartheta_d)\in\Theta$, and $\vj=(j_1,\ldots,j_d)$
where $j_1$,..., $j_d$ are nonnegative integers such that $m:=j_1+\cdots+j_d\leq d+1$, let us define Fourier coefficients of $f$ by
\[\begin{array}{l}
C_{\vk}(x)\;\; :=\frac{1}{(2\pi)^d}\int_{\overline{\cal K}_0}f(x,\vartheta )e^{-i\langle\vk|\vartheta\rangle}d\vartheta,\\
C^{(\vj )}_{\vk}(x) :=\frac{1}{(2\pi)^d}\int_{\overline{\cal K}_0}\frac{\partial^mf}{\partial\vartheta_1^{j_1}\cdots
\partial\vartheta_d^{j_d}}(x,\vartheta )e^{-i\langle\vk |\vartheta\rangle}d\vartheta .
\end{array}\]
Let $\vk^{\vj}:=k_1^{j_1}\cdots k_d^{j_d}$. Since (B ${\cal K}_0$) holds, it is well known 
that
$C^{(\vj )}_{\vk}(x)=i^m \vk^{\vj} C_{\vk}(x)$ for each fixed $x\in E$ (see e.g.\ \cite{Taylor}, pp.\ 177-178). This relation is used in the proof of the next few lemmas (see Appendix).

\begin{lemma}\label{lema4} Let $x,y\in E$. Then for all $\vk\in\Z^d$,
\[\begin{array}{l}
|C_{\vk} (x)|\leq g (x)
\left(\frac{d+1}{1+|k_1|+\cdots+|k_d|}\right)^{d+1}. 
\end{array}\]
\end{lemma}

\begin{lemma}\label{lema:4mom} Let $f\in C(E)$. 
Then for all $0\leq t_0<t$,
\[\begin{array}{ll}
\E{\left(\!\!\int_{t_0}^t\!\! f(X_s)\, dW_s\!\right)\!\!}^4 &\leq 3e^{3(t-t_0)}\E\!\int_{t_0}^t\! f^4 (X_s )\, ds\leq\\
&\leq 24(e^{3(t-t_0)}\E\!\!\int_{t_0}^t\!\! (f(X_s )\! -\! f(X_{t_0}))^4\, ds+\E[f^4(X_{t_0})](t-t_0)^2).
\end{array}\]
\end{lemma}

\begin{lemma}\label{lema1} Let  
{\sc (B4)} hold.  If $c_0 :=(1+c)/2$ then
\begin{equation}\begin{array}{l}
(\forall t\geq 0)\;\; \sup_{0< h\leq h_0}\E\left(\frac{b(X_{t+h})}{b(X_t)}\right)^8 \leq \E\, c_0 (X_t).
\end{array}\label{b1}
\end{equation}
\end{lemma}

\begin{lemma}\label{lema6} There exist constants $K_1>0$, $K_2>0$, $T_0\geq 0$, and $n_0\in\N$, depending on ${\cal K}_0$, $g$ and $a$,
and such that for all $\vk\in\Z^d$, $T>T_0$, $n\geq n_0$ and subdivisions $0=t_0<t_1<\cdots <t_n=T$ (with $\delta_{n,T}\downarrow 0$)
the following hold:
\begin{eqnarray}
\|\!\frac{1}{T\!\sqrt{\delta_{n,T}}}\!\!\sum_{i=0}^{n-1}\!\!\!\int_{t_i}^{t_{i+1}}\!\!\!\!\!\!\!\!\!\!
(C_{\vk}(X_t)\! -\! C_{\vk}(X_{t_i})) a(X_t)\,dt\|_{L^2}\!
&\leq &\! K_1\cdot K_\vk\label{ineq3}\\
\|\!\frac{1}{\sqrt{T\delta_{n,T}}}\!\!\sum_{i=0}^{n-1}\!\!\!\int_{t_i}^{t_{i+1}}\!\!\!\!\!\!\!\!\!\!
(C_{\vk}(X_t)\! -\!C_{\vk}(X_{t_i})) a(X_t)\,dW_t \|_{L^2}\!
&\leq &\!K_1\cdot K_\vk\label{ineq4}\\
\|\!\frac{1}{T\!\sqrt{\delta_{n,T}}}\!\!\sum_{i=0}^{n-1}\!\!\!\int_{t_i}^{t_{i+1}}\!\!\!\!\!\!\!\!\!\!
C_{\vk}(X_{t_i})\!\left(\!\frac{b(X_t )}{b(X_{t_i})}\!-\!1\!\right)\! a(X_t)\,dt\|_{L^2}\!
&\leq &\! K_2\cdot K_\vk\label{ineq5}\\
\|\!\frac{1}{\sqrt{T\delta_{n,T}}}\!\!\sum_{i=0}^{n-1}\!\!\!\int_{t_i}^{t_{i+1}}\!\!\!\!\!\!\!\!\!\!
C_{\vk}(X_{t_i})\!\left(\!\frac{b(X_t )}{b(X_{t_i})}\!-\!1\!\right)\! a(X_t)\,dW_t \|_{L^2}\!
&\leq &\!K_2\cdot K_\vk,\label{ineq6}
\end{eqnarray}
where $K_\vk :=((d+1)/(1+|k_1|+\cdots+|k_d|))^{d+1}$.
\end{lemma}

Let $S_N(x,\vartheta ):=\sum_{|\vk|\leq N} C_{\vk}(x)e^{i\langle\vk |\vartheta\rangle}$ for $x\in E$,
$\vartheta\in\overline{\cal K}_0$ and $N$ be a positive integer.  Then it can be proved that
$\lim_N|S_N(x,\vartheta )-f(x,\vartheta )|=0$ uniformly in $\vartheta\in\overline{\cal K}_0$
by the methods of Fourier analysis (see e.g.\ \cite{Taylor}, pp.\ 180-183).

\begin{lemma}\label{lema5} 
$\sum_{\vk\in\Z^d}|C_{\vk}(x)|\leq K g (x)$, and
$\sup_{N, \vartheta\in\overline{\cal K}_0}|S_N(x,\vartheta )-f(x,\vartheta )|\leq K g (x)$ for a positive and finite constant
\begin{equation}\begin{array}{l}
K=\sum_{\vk\in\Z^d}{\left(\!\!\frac{d+1}{1+|k_1|+\cdots+|k_d|}\!\right)\!.}^{\!\! d+1}
\end{array} \label{K}
\end{equation}
\end{lemma}


\begin{lemma}\label{lema3} Let $a\in C^1 (E)$ and let $f$ be a function that satisfies {\sc (B2)}. Then
for a.s. $\omega\in\Omega$, function $\vartheta\mapsto\int_0^T f(X_t,\vartheta ) a(X_t)\, dW_t(\omega )$ is continuous on $\Theta$.
\end{lemma}

{\em Proof of Theorem \ref{tm:new}}.
Let us prove (\ref{ineq2}) and (\ref{b3}). The proofs of
(\ref{ineq1}) and (\ref{b2}) go in the same way but we have to obtain expressions
of form (\ref{L2}) below with respect to
Lesbegues' instead of Winner's integral, and to apply Lemma \ref{lema6} (\ref{ineq3}) and (\ref{ineq5}).
Without loosing generality let us assume that
${\cal K}={\cal K}_0=\prod_{i=1}^d\langle -\pi,\pi\rangle$ and let $f$ satisfy (B2-3) and (B ${\cal K}_0$).
For fixed $\vartheta\in {\cal K}_0$, $T>0$ and a subdivision
$0=t_0<t_1<\cdots <t_n=T$ we define the following processes:
\[\begin{array}{l}
U_t\quad:=\sum_{i=0}^{n-1}(f(X_t,\vartheta)-f(X_{t_i},\vartheta )) a(X_t)\1_{\langle t_i,t_{i+1}]}(t),\; t\in [0,T],\\
U_t^{(N)}:=\sum_{i=0}^{n-1}(S_N(X_t,\vartheta)-S_N(X_{t_i},\vartheta )) a(X_t)\1_{\langle t_i,t_{i+1}]}(t),\; t\in [0,T],\;\; N\in\N,
\end{array}\]
and
\[\begin{array}{l}
V_t\quad:=\sum_{i=0}^{n-1}f(X_{t_i},\vartheta )\left(\frac{b(X_t)}{b(X_{t_i})}-1\right) a(X_t)\1_{\langle t_i,t_{i+1}]}(t),\; t\in [0,T],\\
V_t^{(N)}:=\sum_{i=0}^{n-1}S_N(X_{t_i},\vartheta )\left(\frac{b(X_t)}{b(X_{t_i})}-1\right) a(X_t)\1_{\langle t_i,t_{i+1}]}(t),\; t\in [0,T],\;\; N\in\N.
\end{array}\]
Then $\lim_N |U_t^{(N)} -U_t|=0$, $\lim_N |V_t^{(N)} -V_t|=0$, and $\sup_N |U_t^{(N)}-U_t|\leq K^2 (g^2 (X_t) + g^2 (X_{t_i})) + a^2 (X_t)/2$, $\sup_N |V_t^{(N)}-V_t|\leq (K/2)g^2 (X_{t_i})\,a^2 (X_t)+2(b(X_t)/b(X_{t_i}))^2+2$,
for $ t\in \langle t_i,t_{i+1}]$ by Lemma \ref{lema5}. Since (B1-4) hold and hence Lemma \ref{lema1} holds there exist $T_1\geq 0$ and $n_1\in\N$ such that for all $T>T_1$, $n\geq n_1$ integrals
$\int_0^T g^2 (X_t)\, dW_t$,  $\sum_{i=0}^{n-1} g^2 (X_{t_i})\Delta_i W$,
$\sum_{i=0}^{n-1} g^2 (X_{t_i})\int_{t_i}^{t_{i+1}}a^2 (X_t)\,dW_t$, $\sum_{i=0}^{n-1}\int_{t_i}^{t_{i+1}}(b(X_t)/b(X_{t_i}))^2\, dW_t$, \\and $\int_0^T a^2 (X_t)\, dW_t$ are well defined, and so
\[\begin{array}{l}
I_{N}(\vartheta ):=\int_0^TU_t^{(N)}\, dW_t\stackrel{\Pb}{\rightarrow}\int_0^T U_t\,dW_t=:I(\vartheta ),\; N\rightarrow +\infty,\\
J_{N}(\vartheta ):=\int_0^TV_t^{(N)}\, dW_t\stackrel{\Pb}{\rightarrow}\int_0^T V_t\,dW_t=:J(\vartheta ),\; N\rightarrow +\infty,
\end{array}\]
by the dominated convergence theorem for stochastic integrals (see e.g.\ \cite{Yor}, Theorem (2.12), pp.\ 134-135).

First, let us consider sequence $(I_N(\vartheta))$. For every
$\vartheta\in {\cal K}_0\cap\Q^d$ 
there exists a subsequence $(N_p)\equiv(N_p(\vartheta))$ and an event $A(\vartheta)$ of the probability 1 such that for all $\omega\in A(\vartheta)$, $\lim_p I_{N_p}(\vartheta)(\omega)=I(\vartheta )(\omega)$. Let us recall that
\[\begin{array}{lcl}
I_{N}(\vartheta)&=&\sum_{i=0}^{n-1}\int_{t_i}^{t_{i+1}}\!(S_N(X_t,\vartheta )-S_N(X_{t_i},\vartheta)) a(X_t)\, dW_t,
\;\; N\in\N,\\
I(\vartheta)&=&
\sum_{i=0}^{n-1}\int_{t_i}^{t_{i+1}}\!(f(X_t,\vartheta )-f(X_{t_i},\vartheta)) a(X_t)\, dW_t.
\end{array}\]
Let $\Omega_0:=\cap_{\vartheta\in {\cal K}_0\cap\Q^d}A(\vartheta )$. Then on this event of probability 1,
for all $\vartheta\in {\cal K}_0\cap\Q^d$, the following holds:
\[\begin{array}{lcl}
|I(\vartheta )| &\leq& |I(\vartheta)-I_{N_p(\vartheta )}(\vartheta )|+|I_{N_p(\vartheta )}(\vartheta )|\leq
|I(\vartheta)-I_{N_p(\vartheta )}(\vartheta )|+\\
&& +\sum_{\vk\in\Z^d}|\sum_{i=0}^{n-1}\int_{t_i}^{t_{i+1}}
(C_{\vk}(X_t)-C_{\vk}(X_{t_i})) a(X_t)\, dW_t|.
\end{array}\]
By taking limit when $p\rightarrow +\infty$, we get the following inequality:
\[\begin{array}{l}
|I(\vartheta )|\leq\sum_{\vk\in\Z^d}|\sum_{i=0}^{n-1}\int_{t_i}^{t_{i+1}}
(C_{\vk}(X_t)-C_{\vk}(X_{t_i}))a(X_t)\, dW_t|.
\end{array}\]
Since $\vartheta\mapsto I(\vartheta )$ is a continuous function by Lemma \ref{lema3}, it turns out that
$\sup_{\vartheta\in\overline{\cal K}_0}|I(\vartheta )|=\sup_{\vartheta\in {\cal K}_0\cap\Q^d}|I(\vartheta )|$, and
so $\sup_{\vartheta\in\overline{\cal K}_0}|I(\vartheta )|$ is a random variable. Hence
\begin{equation}\begin{array}{l}
\sup_{\vartheta\in\overline{\cal K}_0}|I(\vartheta )|\leq\sum_{\vk\in\Z^d}\left|\sum_{i=0}^{n-1}\int_{t_i}^{t_{i+1}}
(C_{\vk}(X_t)-C_{\vk}(X_{t_i})) a(X_t)\, dW_t\right|\;\mbox{\rm  a.s.}
\end{array}\label{L2}
\end{equation}
Since there exist $T_0\geq T_1$ and $n_0\geq n_1$ such that for all $T>T_0$, $n\geq n_0$ and subdivisions of $[0,T]$
with $\delta_{n,T}\downarrow 0$,
\[\begin{array}{l}
\sum_{\vk\in\Z^d}\|\sum_{i=0}^{n-1}\int_{t_i}^{t_{i+1}}(C_{\vk}(X_t)-C_{\vk}(X_{t_i}))a(X_t)\,dW_t\|_{L^2}\leq
K_1 K 
\sqrt{T\delta_{n,T}},
\end{array}\]
by Lemma \ref{lema6} and (\ref{K}), the series on the righthand side of (\ref{L2}) converges a.s.\ and
in $L^2$-norm to a.s.\ equal limits (see Proposition 2.10.1.\ in
\cite{Brockwel}, p.\ 68). Hence $\|\sup_{\vartheta\in {\cal
K}_0}|I(\vartheta )|\|_{L^2}$ $\leq C_1\sqrt{T\delta_{n,T}}$ for $C_1:=K_1 K$. That proves
(\ref{ineq2}). The proof of (\ref{b3}) goes in a similar way considering sequence $(J_N(\vartheta))$. \endproof\vspace{5pt}

We need following lemma for proving consistency and asymptotic normality of
diffusion coefficient parameter estimator.\vspace{5pt}

\begin{lemma}\label{difkoef} Let {\sc (B4)} hold, and let $b\in C^3(E)$. Moreover, let there exist constants $C_b>0$ and $T_b\geq 0$ such that
\[\begin{array}{l}
(\forall T>T_b)\;
\frac{1}{T}\E(\int_0^T ((b^2b''')^2+r^{16})(X_t)\, dt+\sum_{i=0}^{n-1}r^4(X_{t_i})\Delta_it)\leq C_b.
\end{array}\]
Then there exist constants
$C>0$,
$T_0 \geq 0$, and $n_0\in\N$, such that for all $T> T_0$, and $n\geq n_0$,
\[\begin{array}{l}
\frac{1}{T}\E \left|\sum_{i=0}^{n-1}\frac{1}{\Delta_i t}\left(\left(\int_{t_i}^{t_{i+1}}\frac{b(X_t)}{b(X_{t_i})}\, dW_t\right)^2 -(\Delta_i W)^2\right)\right|\leq C.
\end{array}\]
\end{lemma}

\begin{remark}\label{rem4}{\rm If $b$ and its derivatives up to the third order are bounded then the statement of
Lemma \ref{difkoef} hold for all $T>T_0 =0$ by the same arguments as in Remark \ref{rem3}.
}\end {remark}

\subsection{\em Fixed maximal observational time case}

Let $T>0$ be fixed, and let $0=t_0<\cdots t_n=T$, $n\in\N$, be subdivisions of $[0,T]$ such that
$\delta_{n,T}=\max_{0\leq i\leq n-1}\Delta_it\downarrow 0$ when $n\rightarrow +\infty$.
We need the next corollary to Theorem \ref{tm:new}.

\begin{corollary}\label{tm:23}
Let $X$ be a diffusion such that {\sc(H1a-4a)} hold and let
${\cal K}\subset\Theta$ be a relatively compact set. Then for all $\theta_0=(\vartheta_0,\sigma)\in
\Psi$, $T>0$, and $r=0,1,2$,
\begin{equation}\begin{array}{l}
\sup_{\vartheta\in\overline{\cal K}}|D^{r}\ell_{n,T} (\vartheta )-D^{r}\ell_T (\vartheta)|=
O_{\Pb_{\theta_0}}(\sqrt{\delta_{n,T}}),\;\; n\rightarrow +\infty.
\end{array}\label{drugac}
\end{equation}
\end{corollary}

{\em Proof of Corollary \ref{tm:23}}. We prove (\ref{drugac}) for $r=0$.
Statement (\ref{drugac}) for cases $r=1$ and $r=2$ can be proved similarly.
Let  $\theta_0 =(\vartheta_0 ,\sigma)\in\Psi$ be arbitrary, and let
$\mu_0 :=\mu(\cdot,\vartheta_0)$. Moreover, let $f(\cdot,\vartheta):=\mu(\cdot,\vartheta)/b$,
$\vartheta\in\overline{\cal K}$, and $f_0:=\mu_0/b$. Then for any $n$,
\begin{equation}\begin{array}{lcl}
\!\!&\! \!&\! \ell_{n,T} (\vartheta )-\ell_T (\vartheta )  =\\
\!\!&\!=\!&\!\!\sum_{i=0}^{n-1}\!\!\int_{t_i}^{t_{i+1}}\!
 (\frac{\mu (X_{t_i},\vartheta )}{b^2 (X_{t_i})}\! -\!
  \frac{\mu (X_{t},\vartheta )}{b^2 (X_{t})})\, dX_t
-\frac{1}{2}\sum_{i=0}^{n-1}\!\!\int_{t_i}^{t_{i+1}}\!
 (\frac{\mu^2\! (X_{t_i},\vartheta )}{b^2 (X_{t_i})}\! -\!
  \frac{\mu^2\! (X_{t},\vartheta )}{b^2 (X_{t})})\, dt=\\
\!\!&\!=\!&\!\!
 \sum_{i=0}^{n-1}\!\!\int_{t_i}^{t_{i+1}}\!\!( (f(X_t,\vartheta)\!-\!f(X_{t_i},\vartheta))f_0(X_t)
+f(X_{t_i},\vartheta)\left(\!\frac{b(X_t )}{b(X_{t_i})}\!-\!1\!\right)f_0(X_t))\,dt +\\
\!\!&\! \!&\!\! +\sqrt{\sigma}\sum_{i=0}^{n-1}\!\!\int_{t_i}^{t_{i+1}}\!\!( (f(X_t,\vartheta)\!-\!f(X_{t_i},\vartheta))
+f(X_{t_i},\vartheta)\!\left(\!\frac{b(X_t )}{b(X_{t_i})}\!-\!1\!\right)\!)\,dW_t -\\
\!\!&\! \!&\!\! -\frac{1}{2}\sum_{i=0}^{n-1}\!\!\int_{t_i}^{t_{i+1}}\!\!
(f^2(X_t,\vartheta)\!-\!f^2(X_{t_i},\vartheta))\, dt
\end{array}\label{lnld0}
\end{equation}
by the definitions of $\ell_T$ and $\ell_{n,T}$, and (\ref{e1}).

Let us assume for a moment that functions $f_0$, $b$, $b'$, $b''$, are bounded on $E$, and $f$ and
its partial derivatives $D^m_\vartheta f$, $\frac{\partial}{\partial x}D^m_\vartheta f$, and
$\frac{\partial^2}{\partial x^2}D^m_\vartheta f$ are bounded on $E\times\overline{\cal K}$ for $0\leq m\leq d+1$.
Then $f$ and $f^2$
satisfy condition (B2) from Theorem \ref{tm:new}, and  $f_0$ and a constant function $1$ satisfy (B1), since {\sc (H2a-3a)} hold. Hence,
by Remark \ref{rem3} the statements of Theorem \ref{tm:new} holds for these functions, and any $T>0$.
By applying this conclusion to (\ref{lnld0}), the following holds:
\begin{equation}\begin{array}{l}
\|\sup_{\vartheta\in {\cal K}}|\ell_{n,T} (\vartheta )-\ell_T (\vartheta )|\|_{L^2(\Pb_{\theta_0})}
\leq C\sqrt{\delta_{n,T}}, 
\end{array}\label{L2bound}
\end{equation}
for any $T>0$ and subdivisions of $[0,T]$ with $\delta_{n,T}\leq h_0$, and a constant $C>0$ which depends on $T$, $X$ and ${\cal K}$.

Now, let $X$, $\mu$ and $b$ satisfy assumptions {\sc (H1a-3a)}, and let $x_0$ be the initial state of $X$.
Moreover, let $(E_m ,m\geq 1)$ be a sequence
of open and bounded subintervals of $E$ such that for all $m$, $\overline{E}_m\subset E_{m+1}$,
$x_0\in E_1$, and $\bigcup_{m=1}^{+\infty}E_m =E$, and let
$(\phi_m ,m\geq 1)$ be a sequence of
$C^{\infty}$-functions on $E$ such that for all $m$,
$0\leq\phi_m\leq 1$,
$\phi_m (x)=1$ for $x\in\overline{E}_m$ and $\phi_m\equiv0$ on
${E}_{m+1}^c$.  Let us define the following bounded functions for each $m$:
$\mu_m (x,\vartheta ):=
\phi_m (x)\mu (x,\vartheta )$, $(x,\vartheta )\in E\times\Theta$,  $b_m (x):=\phi_m (x)b(x)+c_m (1-\phi_m(x))$,
$x\in E$ where $c_m:=\sign\, b\cdot\max_{x\in\overline{E}_ {m+1}}|b(x)|$.
Since $\mu$ and $b$ satisfy {\sc (Ha2-a3)}, $b_m \in C^2(E)$, and $b_m$, $b_m'$, $b_m''$ are bounded on $E$,
and
$(x,\vartheta)\mapsto\mu_m (x,\vartheta )/b (x)$, $\mu_m^2 (x,\vartheta )/b^2 (x)$  satisfy (B2) and are bounded on $E\times\overline{\cal K}$, and hence satisfy (B3) too, for each $m$.
Moreover, let
$\tau_m :=\inf\{ t\geq 0 : X_t \in E_{m}^c \}$, $m\geq 1$.
Since $X$ is a continuous process, $(\tau_m ,m\geq 1)$
is an increasing sequence of stopping times (see \cite{Yor}) such that
$\tau_m\uparrow +\infty $ a.s., when $m\rightarrow +\infty$.

Let $m$ be fixed and let diffusion $X^m=(X^m_t; t\geq 0)$ be defined as solution to SDE:
\[\begin{array}{l}
X_t^m =x_0 +\int_0^t\mu_m (X_s^m,\vartheta_0 )\, ds+\sqrt{\sigma}\int_0^t b_m (X_s^m )\,
dW_s,\; t>0.
\end{array}\]
By Theorem V.11.2 in \cite{RW} (Vol.\ 2, p.\ 128) such a diffusion exists and is
a.s.\ unique.
Moreover, for almost all $\omega\in\Omega$ and $t\in [0, \tau_m (\omega )]$,
$X_t (\omega )=
X^{m}_t(\omega )$ by Corollary V.11.10 in \cite{RW}  (Vol.\ 2, p.\ 131).
This implies (see \cite{Yor2}) that for an arbitrary number $A>0$,
\begin{equation}\begin{array}{lcl}
&&\!\Pb_{\theta_0}\{\sup_{\vartheta\in {\cal K}}|\ell_{n,T} (\vartheta )-
\ell_T (\vartheta )|>\! A\sqrt{\delta_{n,T}}\}
\leq \\
&\leq&\!\!  \Pb_{\theta_0}\{ \tau_m\!\leq T\}+
\frac{1}{A\sqrt{\delta_{n,T}}}\|\sup_{\vartheta\in {\cal K}}|
\ell_{n,T}^m (\vartheta )-\ell_T^m (\vartheta )|\|_{L^2(\Pb_{\theta_0})},
\end{array}\label{ineqlln}
\end{equation}
where  $\ell_T^m$ and $\ell_{n,T}^m$ are LLF (\ref{ell}) and its Euler approximation (\ref{lnizraz}) respectively, both
based on diffusion $X^m$ with drift $\mu_m (\cdot,\vartheta_0 )$, and diffusion coefficient function
$\sqrt{\sigma}b_m$. Now, (\ref{L2bound}) holds for functions $\ell_T^m$ and $\ell_{n,T}^m$ with constant $C=C_m$. Hence the righthand side
of (\ref{ineqlln}) is dominated by expression
$ \Pb_{\theta_0}\{ \tau_m\!\leq T\}+\frac{1}{A} C_m$.
First, let us take a limit when $n\rightarrow +\infty$, and then when $A\rightarrow +\infty$. Next, we take a limit when $m\rightarrow +\infty$, and hence we prove (\ref{drugac}).
\endproof\vspace{5pt}

{\em Proof of Theorem \ref{tm:30}.\/}
We need to show that the model and random
functions $\ell_T$ and $\ell_{n,T}$, $n\geq 1$, for fixed $T>0$, satisfy conditions (A1-5)
of Theorem 3.1 of \cite{huzakMLE}.
Let ${\cal F}_{n,T}$ be
$\sigma$-subalgebras of ${\cal F}_T^0$ that are introduced in
Section 4. We recall from the same section that $\ell_T$ is a ${\cal
F}_T^0\otimes {\cal B}(\Theta )$-measurable function. In the same
way, $\ell_{n,T}$ is ${\cal F}_{n,T}\otimes {\cal B}(\Theta )$-measurable,
for each $n$. Hence (A1) is satisfied. Corollary  \ref{tm:23} implies
that functions $\ell_T$ and $\ell_{n,T}$, $n\geq 1$, satisfy (A3). The same corollary and {\sc (H5a)} imply (A4) and (A5). Condition (A2) is the same as assumption {\sc (H4a)}. Hence by Theorem 3.1 of \cite{huzakMLE} there exists a sequence of ${\cal F}_T^0$-measurable random vectors
$(\hat{\vartheta}_{n,T} ,n\geq 1)$ such that the statements of Theorem \ref{tm:30} hold.\endproof\vspace{5pt}

For proving Corollary \ref{cor:30a} we need the following lemma.

\begin{lemma}\label{lm:25}
Let {\sc (H1a-2a)} hold, and $T>0$ be fixed. Then for $\theta=(\vartheta,\sigma)\in\Psi$,
\begin{equation}\begin{array}{l}
\sum_{i=0}^{n-1}\frac{(\Delta_i X -\mu (X_{t_i},\vartheta )\Delta_it )^2}{
b^2 (X_{t_i})\Delta_it}-\sigma \sum_{i=0}^{n-1}\frac{(\Delta_i W)^2}{\Delta_it}=O_{\Pb_{\theta}}(1),\;\;n\rightarrow +\infty.
\end{array}\label{n1}
\end{equation}
\end{lemma}

{\em Proof of Corollary \ref{cor:30a}.\/}
Notice that $(ii)$ implies the consistency (i.e.\ $(i)$) of
$\hat{\sigma}_n$. Let us prove $(ii)$. Since
\begin{equation}\begin{array}{l}
\sqrt{n}(\hat{\sigma}_{n,T} -\sigma)=\sqrt{n}(\hat{\sigma}_{n,T} -
\frac{\sigma}{n}\sum_{i=0}^{n-1}\frac{(\Delta_i W)^2}{\Delta_it})+\sigma\sqrt{2}
\cdot\frac{1}{\sqrt{2n}}\sum_{i=0}^{n-1}\frac{(\Delta_i W)^2 -\Delta_it}{\Delta_it}
\end{array}\label{dekomp}
\end{equation}
and $({\cal L})\lim_n
\frac{1}{\sqrt{2n}}\sum_{i=0}^{n-1}\frac{(\Delta_i W)^2 -\Delta_it}{\Delta_it} =N(0,1)$,
 for $(ii)$ to hold it is sufficient to prove 
that for all $\epsilon >0$,
\begin{equation}\begin{array}{l}
\lim_n \Pb_{\theta}\{\sqrt{n}(\hat{\sigma}_{n,T} -
\frac{\sigma}{n}\sum_{i=0}^{n-1}\frac{(\Delta_i W)^2}{\Delta_it})\geq\epsilon\} =0.
\end{array}\label{doh}
\end{equation}
Let $\epsilon >0$ and $\eta >0$ be any numbers and let ${\cal K}$
be a relatively compact set in $\Theta$. If ${\cal K}+\eta:=\{\vartheta\in\Theta : (\exists \vartheta'\in {\cal K})\,
|\vartheta -\vartheta'|<\eta\}$ then
on event
\[\begin{array}{c}
A=
\{|\frac{1}{\sqrt{n}}(\sum_{i=0}^{n-1}\frac{(\Delta_i X -\mu (X_{t_i},
\vartheta)\Delta_it )^2}{b^2 (X_{t_i}) \Delta_it}-\sigma\sum_{i=0}^{n-1}
\frac{(\Delta_i W)^2}{\Delta_it})|<\frac{\epsilon}{5},\; 
\hat{\vartheta}_T\in {\cal K}\}\cap \\ \cap
\{|\hat{\vartheta}_{n,T} -\hat{\vartheta}_T| <\eta,\; 
|\ell_T
(\vartheta)-\ell_T (\hat{\vartheta}_T)|<\sqrt{n}\frac{\epsilon}{10},\; 
|\ell_{n,T} (\vartheta)-\ell_T (\vartheta)|<\sqrt{n}\frac{\epsilon}{10}\}
\cap\\
\cap\{\sup_{\vartheta'\in {\cal K}+\eta }| D\ell_T (\vartheta' )| <
\frac{\sqrt{n}}{\eta}\frac{\epsilon}{10},\; 
\sup_{\vartheta'\in
{\cal K}+\eta }|\ell_{n,T} (\vartheta' )-\ell_T (\vartheta' )|<
\sqrt{n}\frac{\epsilon}{10}\} ,
\end{array}\]
the following holds: $|\sqrt{n}(\hat{\sigma}_{n,T} -\frac{\sigma}{n}\sum_{i=0}^{n-1}\frac{(
\Delta_i W)^2}{\Delta_it})|<\epsilon$.
This implies that
$A\subseteq\{\sqrt{n}(\hat{\sigma}_{n,T}-\frac{\sigma}{n}\sum_{i=0}^{n-1}\frac{(
\Delta_i W)^2}{\Delta_it})<\epsilon\}$.
Hence
\[\begin{array}{rl}
 &\Pb_{\theta}\{\sqrt{n}(\hat{\sigma}_{n,T} -\frac{\sigma }{n}
\sum_{i=0}^{n-1}\frac{(\Delta_i W)^2}{\Delta_it})\geq\epsilon\}\leq\\
\leq &\Pb_{\theta}\{|\frac{1}{\sqrt{n}}(\sum_{i=0}^{n-1}\frac{(\Delta_i X -
\mu (X_{t_i},\vartheta )\Delta_it )^2}{b^2 (X_{t_i}) \Delta_it}-\sigma\sum_{i=0}^{n-1}
\frac{(\Delta_i W)^2}{\Delta_it})|\geq\frac{1}{5}\epsilon\} +\\
 &+\Pb_{\theta}\{\hat{\vartheta}_T\in {\cal K}^c\}+\Pb_{\theta}
\{ |\ell_T (\vartheta )-
\ell_T (\hat{\vartheta}_T)|\geq\sqrt{n}\frac{\epsilon}{10}\}
+\Pb_{\theta}\{ |\hat{\vartheta}_{n,T} -\hat{\vartheta}_T|\geq\eta\}+\\
 &+\Pb_{\theta}\{\sup_{\vartheta'\in {\cal K}+\eta}| D\ell_T (\vartheta' )|
\geq\frac{\sqrt{n}}{\eta}\frac{\epsilon}{10}\} +\Pb_{\theta}\{
|\ell_{n,T}(\vartheta)-\ell_T(\vartheta )| \geq\sqrt{n}\frac{\epsilon}{10}\}
+\\ &+\Pb_{\theta}\{\sup_{\vartheta'\in {\cal K}+\eta}
|\ell_{n,T}(\vartheta' )-\ell_T (\vartheta' )|\geq
\sqrt{n}\frac{\epsilon}{10}\} .
\end{array}\]
By Lemma \ref{lm:25}, Corollary \ref{tm:23}, property $(ii)$ of
$\hat{\vartheta}_{n,T}$ from Theorem \ref{tm:30}, and arbitrariness of ${\cal K}$,  (\ref{doh}) follows.
\endproof 

\subsection{\em Ergodic case}

For all $T>0$ let $0=t_0<\cdots <t_n=T$, $n\in\N$, be equidistant subdivisions of $[0,T]$ such that
$\delta_{n,T}=T/n\rightarrow 0$ when $T\rightarrow +\infty$ and $n\rightarrow +\infty$.
We need the following corollary to Theorem \ref{tm:new}.

\begin{corollary}
\label{tm:cor2}
Let $X$ be a diffusion such that {\sc(H1b-3b)} hold. Then for all $\theta_0=(\vartheta_0,\sigma)\in
\Psi$, $\pi_{\vartheta_0}$-a.s.\ nonrandom initial conditions, and $r=0,1,2$,
\begin{equation}\begin{array}{l}
\sup_{\vartheta\in\overline{\Theta}}|\frac{1}{T}D^{r}\ell_{n,T} (\vartheta )-\frac{1}{T}D^{r}\ell_T (\vartheta)|=
O_{\Pb_{\theta_0}}(\sqrt{\delta_{n,T}}),\;\;T\rightarrow +\infty,\; n\rightarrow +\infty.
\end{array}\label{drugac2}
\end{equation}
\end{corollary}

{\em Proof of Corollary \ref{tm:cor2}}. Similarly to the proof of Corolarlly \ref{tm:23} it is sufficient to prove (\ref{drugac2}) for $r=0$ since the statement of the corollary for cases $r=1$ and $r=2$ can be proved in the same way.
Let  $\theta_0 =(\vartheta_0 ,\sigma)\in\Psi$ be arbitrary, and let
$\mu_0 :=\mu(\cdot,\vartheta_0)$, $\nu:=\sqrt{\sigma} b$, and $\Pb\equiv\Pb_{\theta_0}$, $\E\equiv \E_{\theta_0}$. Let us recall expression (\ref{lnld0}) from the proof of Corolarlly \ref{tm:23} where
$f(\cdot,\vartheta)=\mu(\cdot,\vartheta)/b$,
$\vartheta\in\overline{\Theta}$, and $f_0 =\mu_0/b$. Notice that  $f$ and $f^2$ satisfy (B2) since {\sc (H2a-3a)} hold by {\sc (H2b-3b)}. Let us show
that $f_0$ satisfies (B1) and $f$ satisfies (B3) with respect to $a\equiv f_0$ and compact $\overline{\Theta}$, and that $f^2$ satisfies (B3) with respect to constant function $a\equiv 1$ and the same compact (notice that constant function trivially satisfies (B1)). If we fix $\vartheta\in\overline{\Theta}$, $m$ such that $0\leq m\leq d+1$, and
nonnegative integers $j_1$,..., $j_d$ such that $j_1+\cdots +j_d=m$ then let $\tilde{f}:=\frac{\partial^m}{\partial\vartheta_1^{j_1}\cdots\partial\vartheta_d^{j_d}}f(\cdot,\vartheta )$, and
$\tilde{\mu}:=\frac{\partial^m}{\partial\vartheta_1^{j_1}\cdots\partial\vartheta_d^{j_d}}\mu (\cdot,\vartheta )$. By {\sc (H3a)}, $\tilde{f},\tilde{\mu}\in C^2 (E)$. Since {\sc (H3b)} holds it follows that $|\tilde{f}|\leq g_0\in L^{32}(\pi_{\vartheta_0})\subset
L^{8}(\pi_{\vartheta_0})$, and
\[\begin{array}{lcl}
|\tilde{f}'b|&=&|\tilde{\mu}'-\tilde{f}b'|\leq g_1+g_0 |b'|=:g_{01}\in L^{16}(\pi_{\vartheta_0})\subset L^{8}(\pi_{\vartheta_0})\\
|\tilde{f}'\mu_0|&=&|(\tilde{f}'b)f_0|\leq g_{01}g_1 =:g_{02} \in L^{8}(\pi_{\vartheta_0})\\
|\tilde{f}''b^2|&=&|\tilde{\mu}'' b-2(\tilde{f}'b)b'-f(b''b)|\leq g_2+2g_{01}|b'|+g_0 |b''b|=:g_{03}\in L^{8}(\pi_{\vartheta_0})
\end{array}\]
by {\sc (H2b-3b)}. Then function $g_{00}:=g_0+\sqrt{\sigma}g_{01}+g_{02}+\sigma g_{03}$ is such that  $g_{00}\in L^{8}(\pi_{\vartheta_0})\subset L^{4}(\pi_{\vartheta_0})$
and
\[\begin{array}{l}
\sup_{\vartheta\in\overline{\Theta}}\left(\!|D^m_\vartheta f (\cdot ,\vartheta )|_\infty\!+\!|\frac{\partial}{\partial x}D^m_\vartheta f (\cdot ,\vartheta )|_\infty(|\mu_0 |\! +\!|\nu |)\!+\!
|\frac{\partial^2}{\partial x^2}D^m_\vartheta f (\cdot ,\vartheta  )\nu^2 |_\infty\right)\!\! \leq\! g_{00}
\end{array}\]
for all $0\leq m\leq d+1$. This implies that $f$ satisfies the first part of (B3) with $g\equiv g_{00}$. This also implies that $|f_0|+|\bar{f}_0|\leq g_{00}$ and hence $f_0$, $\bar{f}_0\in L^{8}(\pi_{\vartheta_0})$.
By Chacon-Ornstein theorem, ergodic theorem for additive functionals and its corollary (e.g. Theorem (A.5.2) on p. 504, Theorem (X.3.12)
on p. 397, and Exercise (X.3.18) on p. 399 in \cite{Yor}), for $\pi_{\vartheta_0}$-a.s. initial values $x_0\in E$,
\begin{equation}\begin{array}{lcl}
&&\lim_{T\rightarrow +\infty}\E (\frac{1}{T}\int_0^T f_0^8(X_t)\, dt)=
\lim_{n,T}\E (\frac{1}{T}\sum_{i=0}^{n-1}f_0^8(X_{t_i})\Delta_i t)=\\
&=&\lim_n\E (\frac{1}{n}\sum_{i=0}^{n-1}f_0^8(X_{t_i}))=
\int_E f_0^8 (x)\pi_{\vartheta_0}(dx)\leq \int_E g_{00}^8 (x)\pi_{\vartheta_0}(dx)<+\infty
\end{array}\label{f0}
\end{equation}
since {\sc (H1b)} holds, and subdivisions are equidistant ($\Delta_i t=T/n$ for each $i$). Moreover, since $f_0\in L^{4}(\pi_{\vartheta_0})$ too, the same holds for 4th powers of $f_0$, i.e. if we substitute $f_0^4$ instead of $f_0^8$
in (\ref{f0}). Finally, the both conclusions hold for $\bar{f}_0$ too. Hence $f_0$ satisfies (B1). It remains to show
that $g_{00}$ satisfies the limiting properties from (B3).
Using the same arguments as in proving (\ref{f0}) it follows that (\ref{f0}) holds for 8th and hence for 4th power of
$g_{00}$.
Moreover, since $f_0, g_{00}\in L^{8}(\pi_{\vartheta_0})$ implies $f_0g_{00}\in L^{4}(\pi_{\vartheta_0})$,
and (\ref{f0}) (with respect to $\bar{f}_0$ and $g_{00}$ too) 
holds, it follows that
\[\begin{array}{lcl}
&&
\lim_{n,T}\E (\frac{1}{T}\sum_{i=0}^{n-1}(f_0 g_{00})^4(X_{t_i})\Delta_i t) = 
\int_E (f_0g_{00})^4 (x)\pi_{\vartheta_0}(dx)<+\infty,\\
&&\overline{\lim}_{n,T}\E (\frac{1}{T}\sum_{i=0}^{n-1}g_{00}^4(X_{t_i})\int_{t_i}^{t_{i+1}}(f_0^4 +\bar{f}_0^4)(X_t)\, dt \leq\\
&\leq& \frac{1}{2}\lim_{n,T}\!\E (\frac{1}{T}\!\sum_{i=0}^{n-1}\!\! g_{00}^8(X_{t_i})\Delta_i t)+
\lim_{T\rightarrow +\infty}\!\E (\frac{1}{T}\!\int_0^T\!\!(f_0^8 + \bar{f}_0^8) (X_t)dt)\!<\!+\infty.
\end{array}
\]
Hence $f$ satisfies (B3) for $\pi_{\vartheta_0}$-a.s. nonrandom initial conditions. It remains to show that $f^2$ satisfies
(B3) with respect to function $a\equiv 1$. Let $g:=7\cdot 2^{d+1}g_{00}^2\in L^{4}(\pi_{\vartheta_0})$.
Notice that uniformly with respect to $\vartheta\in\overline{\Theta}$,
\[\begin{array}{l}
|f^2|\!+\!|\frac{\partial}{\partial x}(f^2)|\!+\!|\frac{\partial^2}{\partial x^2}(f^2)|\leq
|f^2|\!+\!2|f\frac{\partial}{\partial x}f|\!+\!2|\left(\frac{\partial}{\partial x}f\right)^2\!+\!f\frac{\partial^2}{\partial x^2}f|\leq
7g_{00}^2\leq g.
\end{array}\]
Let us put $\hat{f}:=\frac{\partial^m}{\partial\vartheta_1^{j_1}\cdots\partial\vartheta_d^{j_d}}(f^2)(\cdot,\vartheta )$ for
fixed $\vartheta\in\overline{\Theta}$, $m$ such that $0\leq m\leq d+1$, and
nonnegative integers $j_1$,..., $j_d$ such that $j_1+\cdots +j_d=m$. Then by induction
\[\begin{array}{l}
|\hat{f}|+|\frac{\partial}{\partial x}\hat{f}|+|\frac{\partial^2}{\partial x^2}\hat{f}|\leq
7\cdot 2^m g_{00}^2\leq g.
\end{array}\]
Then 
(\ref{f0}) (for 4th powers of $g_{00}$) implies that $f^2$ satisfies (B3) with respect to $a\equiv 1$, for $\pi_{\vartheta_0}$-a.s. nonrandom initial conditions. Finally, (B4) holds for $\pi_{\vartheta_0}$-a.s. nonrandom initial conditions since {\sc (H1b-H2b)} hold. Hence we
can apply Theorem \ref{tm:new} to (\ref{lnld0}) to conclude that there exists constants $C>0$, $T_0\geq 0$, and $n_0\in\N$,
such that for all $T>T_0$ and $n\geq n_0$, and arbitrary $A>0$,
\[\begin{array}{cl}
&\Pb_{\theta_0}\{\frac{1}{\sqrt{\delta_{n,T}}}
\sup_{\vartheta\in\overline{\Theta}}|\frac{1}{T}\ell_{n,T}(\vartheta )-\frac{1}{T}\ell_T (\vartheta )|\geq A \}\leq\\
\leq & \frac{1}{A^2}\E{\left(\frac{1}{\sqrt{\delta_{n,T}}}\sup_{\vartheta\in\overline{\Theta}}|\frac{1}{T}\ell_{n,T}(\vartheta )-\frac{1}{T}\ell_T (\vartheta )|\!\right)\!}^2\leq\frac{C}{A^2}.
\end{array}\]
Hence
\[\begin{array}{l}
\lim_{A\rightarrow +\infty}\overline{\lim}_{n,T}\Pb_{\theta_0}\{\frac{1}{\sqrt{\delta_{n,T}}}
\sup_{\vartheta\in\overline{\Theta}}|\frac{1}{T}\ell_{n,T}(\vartheta )-\frac{1}{T}\ell_T (\vartheta )|\geq A \} = 0
\end{array}\]
which proves the corollary. \endproof\vspace{5pt}

In order to prove Theorems \ref{tm:cont}-\ref{tm:erg} we need the following lemmas.

\begin{lemma}\label{lema1b}
Let {\sc (H1b-3b)} hold. Then for all $\theta_0=(\vartheta_0,\sigma)\in\Psi$ there exist constants $C_r>0$ ($r=0,1,2$)
such that $\Pb_{\theta_0}$-a.s. there exists $T_0>0$ such that for all $\vartheta_1,\vartheta_2\in\overline{\Theta}$, and all $T\geq T_0$,
\[\begin{array}{cl}
|\frac{1}{T}D^r\ell_T (\vartheta_1)-\frac{1}{T}D^r\ell_T (\vartheta_2)| &\leq C_r |\vartheta_1-\vartheta_2|,\; r=0,1,2,\\
|\frac{1}{T}D\ell_T (\vartheta_1)-\frac{1}{T}D\ell_T (\vartheta_2)-
\frac{1}{T}D^2\ell_T (\vartheta_2)(\vartheta_1-\vartheta_2)|&\leq\frac{1}{2} C_2 |\vartheta_1-\vartheta_2|^2,\; \mbox{ and}\\
\sup_{\vartheta\in\overline{\Theta}}\frac{1}{T}|D^3\ell_T (\vartheta)| & \leq C_2.
\end{array}\]
\end{lemma}

\begin{lemma}\label{lema2b}
Let {\sc (H1b-3b)} hold. Then for all $\theta_0=(\vartheta_0,\sigma)\in\Psi$, $\Pb_{\theta_0}$-a.s.
\[\begin{array}{l}
\lim_{T\rightarrow +\infty}\sup_{\vartheta\in\overline{\Theta}}|\frac{1}{T}\ell_T (\vartheta )-\ell_{\vartheta_0}(\vartheta )|=0.
\end{array}\]
\end{lemma}

{\em Proof of Theorem \ref{tm:cont}}. Let $\theta_0=(\vartheta_0,\sigma)\in\Psi$ be arbitrary.
Since $\Theta$ is an open set there exists $\varepsilon_0>0$ such that $K(\vartheta_0,\varepsilon_0)\subset\Theta$.
Let $\ell_{\vartheta_0}$ be function (\ref{liml0}) 
and let $
\lambda_0:=\min_{|y|=1}y^\tau I(\vartheta_0 )y=-\max_{|y|=1}y^\tau D^2 \ell_{\vartheta_0}(\vartheta_0 )y >0$
be the minimal eigenvalue of the Fisher information matrix $I(\vartheta_0)$ since it is positive definite by {\sc (H5b)}. Moreover, let $C_r>0$ ($r=0,1,2$) be constants from Lemma \ref{lema1b}, and
let $\Omega_0$ be an intersection of the events from Lemmas \ref{lema1b}-\ref{lema2b}, and the events such that (\ref{liml0}) and
(\ref{liml0dd2}) hold for $\vartheta_0$. Hence $\Pb_{\theta_0} (\Omega_0)=1$, and for $\omega\in\Omega_0$, let $T_0\equiv T_0 (\omega)>0$ be a such that the statements of Lemma \ref{lema1b} hold for $T\geq T_0$.
Let $\varepsilon >0$ be such that $\varepsilon\leq \varepsilon_0\wedge\lambda_0/(4 C_2)$. Then $K(\vartheta_0,\varepsilon)\subset\Theta$. Let $\omega\in\Omega_0$ be fixed.
Since (\ref{liml0dd2}) holds, there exists $T_1\geq T_0$ such that for all
$T\geq T_1$, $|\frac{1}{T}D^2\ell_T (\vartheta_0)-D^2\ell_{\vartheta_0} (\vartheta_0)| <\frac{\lambda_0}{4}$ and
$|\frac{1}{T}D\ell_T (\vartheta_0)-D\ell_{\vartheta_0} (\vartheta_0)| <\frac{\lambda_0}{4}\varepsilon$.
Then for all $y\in\R^d$, $|y|=1$, $T\geq T_1$, and $\vartheta\in K(\vartheta_0,\varepsilon)$,
\[\begin{array}{lcl}
y^\tau (\frac{1}{T}D^2\ell_T (\vartheta )) y &\leq &|\frac{1}{T}D^2\ell_T (\vartheta )-\frac{1}{T}D^2\ell_T (\vartheta_0)|+
|\frac{1}{T}D^2\ell_T (\vartheta_0)-D^2\ell_{\vartheta_0} (\vartheta_0)|+\\
&&+ y^\tau D^2\ell_{\vartheta_0} (\vartheta_0 ) y <C_2 |\vartheta-\vartheta_0|+\frac{\lambda_0}{4} -\lambda_0\leq\\
&\leq & C_2 \frac{\lambda_0}{4C_2}+\frac{\lambda_0}{4} -\lambda_0 =-\frac{\lambda_0}{2}.
\end{array}\]
Hence $\vartheta\mapsto\frac{1}{T}\ell_T(\vartheta)$ is a strictly concave function on $K(\vartheta_0,\varepsilon)$. Moreover,
if $z\in\R^d$ is such that $|z|=\varepsilon$, then for $y:=z/|z|$ and $T\geq T_1$,
\[\begin{array}{lcl}
\frac{1}{T}D\ell_T (\vartheta_0+z ) z &= &\frac{1}{T}D\ell_T (\vartheta_0 )z+z^\tau(\frac{1}{T}\!\!\int_0^1\!\!
D^2\ell_T (\vartheta_0+tz)\,dt)z\leq\\
&\leq&|\frac{1}{T}D\ell_T (\vartheta_0 )-D\ell_{\vartheta_0}(\vartheta_0)|\varepsilon +y^\tau(\frac{1}{T}\!\!\int_0^1\!\!
D^2\ell_T (\vartheta_0+tz)dt)y\varepsilon^2\leq\\
&\leq&\frac{\lambda_0}{4}\varepsilon^2 -\frac{\lambda_0}{2}\varepsilon^2 =-\frac{\lambda_0}{4}\varepsilon^2<0.
\end{array}\]
Then there exists $\hat{\vartheta}_T\in K(\vartheta_0,\varepsilon)$ such that $D\ell_T(\hat{\vartheta}_T)=\O$ (see e.g.
Lemma 4.3. in \cite{huzakMLE}), and $D^2\ell_T(\hat{\vartheta}_T)<\O$ since
$\min_{|y|=1}y^\tau (-\frac{1}{T}D^2\ell_T (\vartheta )) y\geq\frac{\lambda_0}{2}=\frac{1}{2}\min_{|y|=1}y^\tau I(\vartheta_0 )y$
for all $\vartheta\in K(\vartheta_0,\varepsilon)$ obviously. Since $\varepsilon >0$ is an arbitrary small number, these imply statement $(ii)$ of the theorem. Notice that $\hat{\vartheta}_T$ is the unique point of maximum of function $\ell_T$ on
$K(\vartheta_0,\varepsilon)$ since $\ell_T$ is strictly concave on this set. To finish the proof of statement $(i)$ we have to prove that there exists $T_2\geq T_1$ such that $\hat{\vartheta}_T$ is the unique point of global maximum of $\ell_T$ on $\overline{\Theta}$. Since for all $\vartheta\in\overline{\Theta}\setminus\{\vartheta_0\}$, $\ell_{\vartheta_0}(\vartheta_0 )>\ell_{\vartheta_0}(\vartheta )$, $\ell_{\vartheta_0} \in C(\overline{\Theta})$, and $\overline{\Theta}\setminus K(\vartheta_0,\varepsilon)$ is a compact set,
it follows that $\ell_{\vartheta_0}(\vartheta_0 )>\sup_{|y|\geq\varepsilon}\ell_{\vartheta_0}(\vartheta_0+y )$. By Lemma 4.4. in \cite{huzakMLE} there exists a number $0<s(\varepsilon)<\varepsilon$ such that
\[\begin{array}{l}
\Delta (\vartheta_0,\varepsilon):=\inf_{|x|\leq s(\varepsilon)}\ell_{\vartheta_0}(\vartheta_0+x )-\sup_{|y|\geq\varepsilon}\ell_{\vartheta_0}(\vartheta_0+y )>0.
\end{array}\]
Since Lemma \ref{lema2b} holds there exists $T_2\geq T_1$ such that for
$T\geq T_2$,
\[\begin{array}{l}
\sup_{\vartheta\in\overline{\Theta}}|\frac{1}{T}\ell_T(\vartheta )-\ell_{\vartheta_0}(\vartheta )|<\frac{\Delta (\vartheta_0,\varepsilon)}{4}.
\end{array}\]
If $x,y\in\R^d$ such that $|x|\leq s(\varepsilon )$ and $|y|\geq\varepsilon $ then
\[\begin{array}{lcl}
\frac{1}{T}\ell_T(\vartheta_0 +x) &=&\frac{1}{T}\ell_T(\vartheta_0 +x)-\ell_{\vartheta_0}(\vartheta_0 +x) +\ell_{\vartheta_0}(\vartheta_0 +x) -\ell_{\vartheta_0}(\vartheta_0 +y)+\\
&&+\ell_{\vartheta_0}(\vartheta_0 +y) -\frac{1}{T}\ell_T(\vartheta_0 +y)+
\frac{1}{T}\ell_T(\vartheta_0 +y)\geq\\
&\geq &-\frac{\Delta (\vartheta_0,\varepsilon)}{4}+\inf_{|x|\leq s(\varepsilon)}\ell_{\vartheta_0}(\vartheta_0+x )-\sup_{|y|\geq\varepsilon}\ell_{\vartheta_0}(\vartheta_0+y )+\\
&& -\frac{\Delta (\vartheta_0,\varepsilon)}{4}+\frac{1}{T}\ell_T(\vartheta_0 +y)\geq\\
&\geq &\frac{\Delta (\vartheta_0,\varepsilon)}{2}+\frac{1}{T}\ell_T(\vartheta_0 +y)
\end{array}\]
implying that
\begin{equation}\begin{array}{l}
\inf_{|x|\leq s(\varepsilon)}\frac{1}{T}\ell_{T}(\vartheta_0+x )-\sup_{|y|\geq\varepsilon}\frac{1}{T}\ell_{T}(\vartheta_0+y )\geq \frac{\Delta (\vartheta_0,\varepsilon)}{2} >0
\end{array}\label{deltapotez}\end{equation}
and hence $\ell_T(\vartheta_0)> \sup_{|y|\geq\varepsilon}\ell_T(\vartheta_0+y )$. Finally, $(i)$ follows. To prove statement $(iii)$,
first notice that
\begin{equation}\begin{array}{l}
\frac{1}{\sqrt{T}}D\ell_{T}(\vartheta_0)=\frac{\sqrt{\sigma}}{\sqrt{T}}\int_0^T\frac{1}{b(X_t)}D\mu_0 (X_t)\, dW_t\stackrel{{\cal L}-\Pb_{\theta_0}}{\longrightarrow}N(\O,\sigma I(\vartheta_0)),\; T\rightarrow +\infty
\end{array}\label{elnorm}\end{equation}
by Theorem 1 in \cite{BH} since {\sc (H1b-5b)} hold, and second notice that for $\bar{\vartheta} (s):=s\vartheta_T +(1-s)\vartheta_0$,
\begin{equation}\begin{array}{l}
D\ell_{T}(\hat{\vartheta}_T)\! =\! D\ell_{T}(\vartheta_0)\!+\!D^2\ell_{T}(\vartheta_0) (\hat{\vartheta}_T\!-\!\vartheta_0)\!+\! \int_0^1\!\!\int_0^1\!
D^3\ell_{T}(\bar{\vartheta}(st))\, ds\, tdt (\hat{\vartheta}_T\!-\!\vartheta_0)^2.
\end{array}\label{eltaylor}\end{equation}
Let $H_T(\vartheta_0):=\frac{1}{T}D^2\ell_{T}(\vartheta_0)+\frac{1}{T}\int_0^1\!\!\int_0^1\!
D^3\ell_{T}(\bar{\vartheta}(st))\, ds\,tdt (\hat{\vartheta}_T\!-\!\vartheta_0)$, and let us recall $\omega\in\Omega_0$ and $T_1=T_1(\omega )$ from the first part od the proof. Notice that $H_T(\vartheta_0)$ is a symmetric matrix. Then from Lemma \ref{lema1b}, for $T\geq T_1$,
\[\begin{array}{l}
|H_T(\vartheta_0)-\frac{1}{T}D^2\ell_T (\vartheta_0)|\leq\sup_{\vartheta\in\overline{\Theta}}|\frac{1}{2T}D^3\ell_{T}(\vartheta )|
|\hat{\vartheta}_T\!-\!\vartheta_0|\leq\frac{C_2}{2}|\hat{\vartheta}_T\!-\!\vartheta_0|
\end{array}\]
and hence, for $y\in\R^d$ such that $|y|=1$,
\[\begin{array}{l}
y^\tau H_T (\vartheta_0 ) y \leq |H_T(\vartheta_0)-\frac{1}{T}D^2\ell_T (\vartheta_0)|+y^\tau (\frac{1}{T}D^2\ell_T (\vartheta )) y\leq 
-\frac{3\lambda_0}{8}
\end{array}\]
implying that $H_T (\vartheta_0 )$ is a negative definite matrix, and $|H_T (\vartheta_0 )^{-1}|\leq\frac{8}{3\lambda_0}$. Since
$|I(\vartheta_0)^{-1}|=1/\lambda_0$,
\[\begin{array}{cl}
&|H_T(\vartheta_0)^{-1}+I(\vartheta_0)^{-1}|\leq |H_T(\vartheta_0)^{-1}|\cdot |H_T(\vartheta_0)+I(\vartheta_0)|
\cdot |I(\vartheta_0)^{-1}| \leq\\
\leq & \frac{8}{3\lambda_0^2}(\frac{C_2}{2}|\hat{\vartheta}_T\!-\!\vartheta_0|+|\frac{1}{T}D^2\ell_T (\vartheta_0)-D^2\ell_{\vartheta_0}(\vartheta_0)|),
\end{array}\]
and $(ii)$ and (\ref{liml0dd2}) hold, it follows that $\Pb_{\theta_0}$-a.s.
\begin{equation}\begin{array}{l}
\lim_{T\rightarrow +\infty}H_T(\vartheta_0)^{-1}=-I(\vartheta_0)^{-1}.
\end{array}\label{limH}\end{equation}
Finally, since $D\ell_{T}(\hat{\vartheta}_T)=\O$ and $I(\vartheta_0)$ is nonrandom, (\ref{elnorm}-\ref{limH}) imply that
\[\begin{array}{l}
\sqrt{T}(\hat{\vartheta}_T-\vartheta_0)=-H_T(\vartheta_0)^{-1}\frac{1}{\sqrt{T}}D\ell_T(\vartheta_0)\stackrel{{\cal L}-\Pb_{\vartheta_0}}{\longrightarrow}
N(\O,\sigma I(\vartheta_0 )^{-1}),\; T\rightarrow +\infty.\;\mbox{\endproof}
\end{array}\]

{\em Proof of Theorem \ref{tm:erg}}. Let $\theta_0=(\vartheta_0,\sigma)\in\Psi$ be arbitrary, and let
$C_r>0$ ($r=0,1,2$) be constants from Lemma \ref{lema1b}. Moreover, let $\Omega_0$ be a $\Pb_{\theta_0}$-probability one
event from Lemmas \ref{lema1b}-\ref{lema2b} and Theorem \ref{tm:cont} ($i$-$ii$). Let $\omega\in\Omega_0$ be
fixed. Let  $\varepsilon_0>0$ be a such number that $K(\vartheta_0,\varepsilon_0)\subset\Theta$, and let $\lambda_0>0$
be the minimal eigenvalue of Fisher matrix $I(\vartheta_0)$.
Then there exists $T_0=T_0(\omega)\geq 0$ such that for all $T>T_0$, $\hat{\vartheta}_T\in K(\vartheta_0,\varepsilon_0/2)$ and $\lambda_T:=\min_{|y|=1}y^\tau (-\frac{1}{T}D^2\ell_T(\hat{\vartheta}_T))y\geq\lambda_0/2>0$, and the statements of Lemma \ref {lema1b} hold.
Let $\varepsilon >0$ be an arbitrary small number such that $\varepsilon<
\frac{\varepsilon_0}{2}\wedge \frac{\lambda_0}{8C_2}$. Then $K(\hat{\vartheta}_T,\varepsilon)\subset K(\hat{\vartheta}_T,\varepsilon_0/2)\subset
K(\vartheta_0,\varepsilon_0)\subset\Theta$. Moreover, on event
\[\begin{array}{l}
\Omega_{n,T}:=\{ \sup_{\vartheta\in\overline{\Theta}}|\frac{1}{T}D^r\ell_{n,T}(\vartheta)-\frac{1}{T}D^r\ell_{T}(\vartheta)|\leq
\frac{\lambda_0}{8}(1\wedge\frac{\lambda_0}{8C_2}),\; r=1,2\},
\end{array}\]
for $\vartheta\in K(\hat{\vartheta}_T,\varepsilon)$ and $z\in\R^d$ such that $|z|=\varepsilon$, and $y:=z/|z|$,
the following holds:
\[\begin{array}{lcl}
y^\tau D^2\ell_{n,T} (\vartheta ) y &\leq &|D^2\ell_{n,T} (\vartheta )-D^2\ell_T (\vartheta)|\!+\!
|D^2\ell_T (\vartheta)-D^2\ell_{T} (\hat{\vartheta}_T)|\!+\\
&&+ y^\tau D^2\ell_{T} (\hat{\vartheta}_T) y <(\frac{\lambda_0}{4}+C_2 \frac{\lambda_0}{8C_2} -\frac{\lambda_0}{2})T =-\frac{\lambda_0}{4}T<0,\\
D\ell_{n,T} (\hat{\vartheta}_T+z ) z &= &D\ell_{n,T} (\hat{\vartheta}_T )z+z^\tau(\!\!\int_0^1\!\!
D^2\ell_{n,T}(\hat{\vartheta}_T+tz)\,dt) z\leq\\
&\leq&|D\ell_{n,T} (\hat{\vartheta}_T )-D\ell_T (\hat{\vartheta}_T )|\varepsilon +y^\tau(\int_0^1\!\!
D^2\ell_T (\hat{\vartheta}_T+tz)dt)y\varepsilon^2\leq\\
&\leq&\varepsilon\frac{\lambda_0}{8C_2}(\frac{\lambda_0}{8} -\frac{\lambda_0}{4})T =-\varepsilon\frac{\lambda_0}{8C_2}\frac{\lambda_0}{8}T<0.
\end{array}\]
Hence $\vartheta\mapsto\ell_{n,T}(\vartheta)$ is a strictly concave function on $K(\hat{\vartheta}_T,\varepsilon)$,
and there exists $\hat{\vartheta}_{n,T}\in K(\hat{\vartheta}_T,\varepsilon)$ such that $D\ell_{n,T}(\hat{\vartheta}_{n,T})=\O$, and $\hat{\vartheta}_{n,T}$ is the unique stationary point and a point of
maximum of $\ell_{n,T}$ at $K(\hat{\vartheta}_T,\varepsilon)$. These imply that $\hat{\vartheta}_{n,T}$ is a random vector. Since $\lim_{n,T}\Pb_{\theta_0}(\Omega_{n,T}^c)=0$ by Corollary \ref{tm:cor2},
and $\Omega_{n,T}\subset\{ D\ell_{n,T}(\hat{\vartheta}_{n,T})=\O\}\cap\{|\hat{\vartheta}_{n,T}-\hat{\vartheta}_T|<\varepsilon\}$, statements $(i)$ and $(ii)$ of the theorem follow. Moreover if process $(\tilde{\vartheta}_{n,T})$ satisfies $(i)$ and $(ii)$ then
statement $(iv)$ follows since
\[\begin{array}{l}
\Omega_{n,T}\cap\{D\ell_{n,T}(\tilde{\vartheta}_{n,T})=\O\}\cap\{|\tilde{\vartheta}_{n,T}-\hat{\vartheta}_T|<\varepsilon\}
\subseteq\{\hat{\vartheta}_{n,T}=\tilde{\vartheta}_{n,T}\}
\end{array}\]
by uniqness of a stationary point of $\ell_{n,T}$ on $K(\hat{\vartheta}_{T},\varepsilon)$. To prove $(iii)$, let $A>0$
be an arbitrary number, and let
$\Omega_{n,T}(A):=\{
\sup_{\vartheta\in\overline{\Theta}}|\frac{1}{T}D\ell_{n,T}(\vartheta)-\frac{1}{T}D\ell_{T}(\vartheta)|\leq
\frac{\lambda_0}{4}A\sqrt{\delta_{n,T}}\}$. Then on event $\Omega_{n,T}(A)\cap\Omega_{n,T}$,
\[\begin{array}{ccl}
|\hat{\vartheta}_{n,T}\!-\!\hat{\vartheta}_T|&\!\leq\! &|(D^2\ell_{T}(\hat{\vartheta}_T))^{-1}|\!\cdot \! |D^2\ell_{T}(\hat{\vartheta}_T)
(\hat{\vartheta}_{n,T}-\hat{\vartheta}_T)|\leq\\
&\!\leq\! &|(D^2\ell_{T}(\hat{\vartheta}_T))^{-1}|\!\cdot \!
|D\ell_{T}(\hat{\vartheta}_{n,T})\!-\!D\ell_{T}(\hat{\vartheta}_T)\!-\!D^2\ell_{T}(\hat{\vartheta}_T)
(\hat{\vartheta}_{n,T}\!-\!\hat{\vartheta}_T)|+\\
&&+|(D^2\ell_{T}(\hat{\vartheta}_T))^{-1}|\!\cdot \!|D\ell_{n,T}(\hat{\vartheta}_{n,T})-D\ell_{T}(\hat{\vartheta}_{n,T})|\leq\\
&\leq&\frac{2}{\lambda_0 T}\frac{C_2}{2}\frac{\lambda_0}{2C_2}T|\hat{\vartheta}_{n,T}-\hat{\vartheta}_T|+\frac{2}{\lambda_0 T}\frac{\lambda_0 T}{4}A\sqrt{\delta_{n,T}}\leq\\
&\leq &\frac{1}{2}|\hat{\vartheta}_{n,T}-\hat{\vartheta}_T|+\frac{1}{2}A\sqrt{\delta_{n,T}}\\
\Rightarrow & & |\hat{\vartheta}_{n,T}\!-\!\hat{\vartheta}_T|\leq A\sqrt{\delta_{n,T}}
\end{array}\]
by Lemma \ref{lema1b} and since $\varepsilon\leq\frac{\lambda_0}{2C_2}$. Hence $\Omega_{n,T}(A)\cap\Omega_{n,T}\subseteq\{
|\hat{\vartheta}_{n,T}\!-\!\hat{\vartheta}_T|\leq A\sqrt{\delta_{n,T}}\}$, and
\[\begin{array}{cl}
0&\leq\overline{\lim}_{A\rightarrow +\infty}\overline{\lim}_{n,T}\Pb_{\theta_0}\{|\hat{\vartheta}_{n,T}\!-\!\hat{\vartheta}_T|\leq A\sqrt{\delta_{n,T}}\}\leq\\
&\leq \lim_{A\rightarrow +\infty}\overline{\lim}_{n,T}\Pb_{\theta_0}(\Omega_{n,T}(A)^c)+\overline{\lim}_{n,T}\Pb_{\theta_0}(\Omega_{n,T}^c)=0
\end{array}\]
by Corollary \ref{tm:cor2}, and $(iii)$ follows. Consistency of $\hat{\vartheta}_{n,T}$ (the first part of statement $(v)$) follows directly from $(ii)$ and Theorem \ref{tm:cont} $(ii)$. To prove its asymptotic normality (the second part of $(v)$) notice that
\[\begin{array}{l}
|\sqrt{T}(\hat{\vartheta}_{n,T}\!-\!\vartheta_0)-\sqrt{T}(\hat{\vartheta}_{T}\!-\!\vartheta_0)|=
\sqrt{T\delta_{n,T}}\frac{1}{\sqrt{\delta_{n,T}}}|\hat{\vartheta}_{n,T}\!-\!\hat{\vartheta}_T|\stackrel{\Pb_{\theta_0}}{
\longrightarrow} 0
\end{array}\]
when $\lim_{n,T}T\delta_{n,T}=0$ since $(iii)$ holds. Then the second part of $(v)$ follows by Slutsky theorem since Theorem \ref{tm:cont} $(iii)$ holds. To prove statement $(vi)$, first we need to prove that
\begin{equation}\begin{array}{l}
\frac{1}{T}(\sum_{i=0}^{n-1}\!\frac{(\Delta_i X -\mu (X_{t_i},\vartheta )\Delta_it )^2}{
b^2 (X_{t_i})\Delta_it}-\sigma\! \sum_{i=0}^{n-1}\!\frac{(\Delta_i W)^2}{\Delta_it})=O_{\Pb_{\theta}}(1),\;T\rightarrow +\infty, n\rightarrow +\infty
\end{array}\label{opb}\end{equation}
for $\pi_{\vartheta_0}$-a.s. initial conditions.
This follows from Lemma \ref{difkoef}, the proof of Lemma \ref{lm:25}, and the fact that the functions $f:=\mu(\cdot,\vartheta)/b$ and $b$ satisfies (B1-4) which is proved in Corollary \ref{tm:cor2}. The proof of asymptotic normality of $\hat{\sigma}_{n,T}$ is the same as in the proof of Corollary \ref{cor:30a} since
\[\begin{array}{cl}
&\frac{1}{\sqrt{n}}(\sum_{i=0}^{n-1}\!\frac{(\Delta_i X -\mu (X_{t_i},\vartheta )\Delta_it )^2}{
b^2 (X_{t_i})\Delta_it}-\sigma\! \sum_{i=0}^{n-1}\!\frac{(\Delta_i W)^2}{\Delta_it})=\\
=&
\frac{\sqrt{T\delta_{n,T}}}{T}(\sum_{i=0}^{n-1}\!\frac{(\Delta_i X -\mu (X_{t_i},\vartheta )\Delta_it )^2}{
b^2 (X_{t_i})\Delta_it}-\sigma\! \sum_{i=0}^{n-1}\!\frac{(\Delta_i W)^2}{\Delta_it})\rightarrow 0
\end{array}\]
when $T\rightarrow +\infty$ such that $T\delta_{n,T}\rightarrow 0$, and since $(i-v)$, Corollary \ref{tm:cor2},
and Lemma \ref{lema2b} hold.
Similarly  consistency of $\hat{\sigma}_{n,T}$
follows from decomposition (\ref{dekomp}) in the proof of Corollary \ref{cor:30a} (but without factor "$\sqrt{n}$") by
using (\ref{opb}) which appears with factor "$\delta_{n,T}$" (notice that $\delta_{n,T}/T =1/n$), and by the strong low of large numbers instead of CLT. In this case it is sufficient to assume that $\delta_{n,T}\rightarrow 0$ when $T\rightarrow +\infty$.
Finally, for proving ${\cal F}_{n,T}^0$-measurability of $\hat{\vartheta}_{n,T}$ (and hence $\hat{\sigma}_{n,T}$ too) it is sufficient to
prove that $\hat{\vartheta}_{n,T}$ is a unique point of maximum of $\ell_{n,T}$ on $\overline{\Theta}$. This proof
follows in the similar way as proof of uniqness of $\hat{\vartheta}_T$ as global point of maximum of $\ell_T$ on $\overline{\Theta}$ by replacing $\ell_T$ with $\ell_{n,T}$ and $\ell_{\vartheta_0}$ with $\ell_T$.\endproof


\begin{center}
{\footnotesize APPENDIX}
\end{center}
\addvspace{.15in}\nopagebreak
   \markboth{APPENDIX}{APPENDIX}
\markright{APPENDIX}

{\footnotesize

{\em Proof of Lemma \ref{lema4}}. 
Let ${\vk}^{\vj}=k_1^{j_1}\cdots k_d^{j_d}$ for nonnegative integers $j_1$,..., $j_d$  such that
$m:=j_1+\cdots +j_d\leq d+1$. Then for $x\in E$,
\[\begin{array}{l}
{|\vk|}^{\vj}
|C_{\vk}(x)|=|C_{\vk}^{(\vj)}(x)|\leq
\frac{1}{(2\pi)^d}\int_{\overline{\cal K}_0}\left|\frac{\partial^m}{\partial\vartheta_1^{j_1}\cdots
\partial\vartheta_d^{j_d}}f(x,\vartheta)\right|d\vartheta \leq g(x)
\end{array}\]
by the definition of Fourier coefficients, the monotonicity of integral, and (B3).
Hence
\[\begin{array}{lcl}
&&(1+|k_1|+\cdots +|k_d|)^{d+1} |C_{\vk}(x)|=\\
&=&\sum_{j_0+j_1+\cdots +j_d= d+1}\frac{(d+1)!}{j_0!j_1!\cdots j_d!}{|\vk|}^{\vj}
|C_{\vk}(x)|\leq (d+1)^{d+1} g(x)
\end{array}\]
by multinomial theorem, which implies the statement of the lemma.
\endproof\vspace{5pt}

{\em Proof of Lemma \ref{lema:4mom}}. At first, let us suppose that $f$ is bounded on $E$. If
$M_t :=(\int_{t_0}^t f(X_s)\, dW_s)^2 -\int_{t_0}^t f^2(X_s)\, ds$, then It\^{o} formula and isometry
implies
\[\begin{array}{l}
\E(M_t)^2=4\E(\int_{t_0}^t (\int_{t_0}^s f(X_u)\,dW_u)f(X_s)\, dW_s)^2\leq 2\|f\|_\infty^4 (t-t_0)^2.
\end{array}\]
Hence, if $N_t:=\int_{t_0}^{t} f(X_s)\,dW_s$ then
\[\begin{array}{l}
\E(N_t)^4\leq 2\E(M_t)^2+2\E(\int_{t_0}^t f^2(X_s)\, ds)^2\leq 6\|f\|_\infty^4 (t-t_0)^2<+\infty.
\end{array}\]
Similarly, if  $t_0\leq s<s+h\leq t$ then
\[\begin{array}{l}
\E(N_{s+h}^2-N_s^2)^2\leq 2\|f\|_\infty^4 (4(s-t_0)+3h) h\rightarrow 0,\;\; h\rightarrow 0.
\end{array}\]
In addition $\E(f^2(X_{s+h})-f^2(X_{s}))^2\rightarrow 0$, and $\E f^4(X_{s+h})\rightarrow \E f^4(X_{s})$ when $h\rightarrow 0$ by
the dominated convergence theorem. Hence
\[\begin{array}{lcl}
&& |\E(N_{s+h}^2f^2(X_{s+h}))-\E(N_s^2f^2 (X_s))|\leq\\
&\leq& \|f\|_\infty^2\E|N_{s+h}^2-N_s^2|+\E(N_s^2|f^2(X_{s+h})-f^2(X_{s})|)\leq\\
&\leq&\|f\|_\infty^2\sqrt{\E(N_{s+h}^2-N_s^2)^2}+\sqrt{\E(N_s^4)\E(f^2(X_{s+h})-f^2(X_{s}))^2}\rightarrow 0,\;\; h\rightarrow 0
\end{array}\]
implying that  $s\mapsto \E(N_s^2f^2 (X_s))$, and $s\mapsto \E f^4(X_s)$ are continuous functions on $[t_0,t]$.
Let $x(t):=\E N_t^4$. Since for $t_0\leq s<s+h\leq t$,
\[\begin{array}{l}
x(s+h)-x(s)=\E(N_{s+h}^4-N_s^4)=6\int_s^{s+h} \E(N_u^2f^2 (X_u))\, du,
\end{array}\]
by It\^{o} formula, $s\mapsto x(s)$ is a differentiable function for $s>t_0$, and
\[\begin{array}{l}
x(s+h)-x(s)\leq 3\int_s^{s+h} x(u)\, du +3\int_s^{s+h}\E f^4(X_u)\, du.
\end{array}\]
Hence $\dot{x}(s)\leq 3x(s)+3\E f^4(X_s)$ for $s>t_0$, and $x(t_0)=0$, implying that
\[\begin{array}{l}
\E(\int_{t_0}^{t}\!\!\! f(X_s)\,dW_s)^4=x(t)\leq 3e^{3t}\!\!\int_{t_0}^t\!\! e^{-3s}\E f^4(X_s)\, ds\leq 3e^{3(t-t_0)}\E(\int_{t_0}^t\!\!\! f^4(X_s)\, ds).
\end{array}\]
Now, let $f\in C(E)$ be unbounded generally. Then there exists a sequence $(f_m)$ of bounded functions
such that for all $m$, $f_m\in C(E)$, $|f_m|\uparrow |f|$, and $f_m \rightarrow f$ (see the proof of Corollary \ref{tm:23}).
Since $\lim_m f_m(X_s)=f(X_s)$ and $|f_m(X_s)-f(X_s)|\leq 2|f(X_s)|$ for all $m$, and $s$, it follows that
\[\begin{array}{l}
\int_{t_0}^t \!\! f_m(X_s)\,dW_s\stackrel{\Pb}{\longrightarrow}\int_{t_0}^t \!\! f(X_s)\,dW_s
\end{array}\]
by the dominated convergence theorem for stochastic integrals. Then there exists a subsequence such that
\[\begin{array}{l}
\int_{t_0}^t \!\! f_{m_k}(X_s)\,dW_s\stackrel{\mbox{\rm a.s.}}{\longrightarrow}\int_{t_0}^t \!\! f(X_s)\,dW_s\;\Rightarrow
{\left(\!\int_{t_0}^t \!\!f_{m_k}(X_s)\,dW_s\right)\!\!}^4\stackrel{\mbox{\rm a.s.}}{\longrightarrow}{\left(\!\int_{t_0}^t \!\!  f(X_s)\,dW_s\!\right)\!\!}^4,
\end{array}\]
and hence by Fatou's lemma and monotone convergence theorem
\[\begin{array}{lcl}
&&\E{\left(\int_{t_0}^t \!\! f(X_s)\,dW_s\right)\!\!}^4\leq \underline{\lim}_k\E{\left(\int_{t_0}^t \!\! f_{m_k}(X_s)\,dW_s\right)\!\!}^4\leq\\
&\leq& 3e^{3(t-t_0)}\overline{\lim}_k\E(\int_{t_0}^t\!\! f^4_{m_k}(X_s)\, ds)
)= 3e^{3(t-t_0)}\E(\int_{t_0}^t\!\! f^4(X_s)\, ds). 
\end{array}\]
The last inequality follows trivially from the first one. \endproof\vspace{5pt}

{\em Proof of Lemma \ref{lema1}}. By applying It\^{o} formula on log-function of $b^{16}$ over time interval $[t,t+h]$ it
follows that $(b(X_{t+h})/b(X_t))^{16} =M_h Z_h$ where
\[\begin{array}{l}
M_h=\exp\left(16\sqrt{\sigma}\int_t^{t+h}b'(X_s)\, dW_s-\frac{\sigma}{2}16^2\int_t^{t+h}b'^2(X_s)\,ds\right)
\end{array}\]
is a positive supermartingal (see \cite{Shiryayev}, Lemma 6.1, p.207) since $\E\int_t^{t+h}\!b'^2(X_s)\, ds\leq\E\int_t^{t+h}\!r^2(X_s)\, ds\leq (3T+\E\int_0^{T}\!r^8(X_s)\,ds)/4<+\infty$
by assumption (B4) (and so $\E M_h\leq \E M_0=1$), and
\[\begin{array}{l}
Z_h =\exp{\left( 8\int_t^{t+h}\left(2\frac{\mu_0 b'}{b} +\sigma (b''b+15b'^2)\right)(X_s)\,ds\right)\!}.
\end{array}\]
By Markov property and assumption (B4), for $0<h\leq h_0$,
\[\begin{array}{cl}
&\E Z_h =\E [\, \E[Z_h | {\cal F}_t^0]\,]=\\
=&\E[\, \E_{X_t}[\,\exp{\left( 8\int_0^{h}\left(2\frac{\mu_0 b'}{b} +\sigma (b''b+15b'^2)\right)(X_s)\,ds\right)}]\,]\leq \E [c(X_t)].
\end{array}\]
Hence, for $0<h\leq h_0$,
\[\begin{array}{l}
\E{\left(\!\frac{b(X_{t+h})}{b(X_t)}\!\right)\!}^8=\E\sqrt{M_h Z_h}\leq\frac{1}{2}(\E M_h+\E Z_h )\leq\frac{1}{2}(1+\E\, c(X_t ))=\E\, c_0(X_t). \mbox{\endproof}
\end{array}\]

{\em Proof of Lemma \ref{lema6}}. First, let us show that (\ref{ineq4}) implies (\ref{ineq3}). In the same way it can be shown that (\ref{ineq6}) implies (\ref{ineq5}). Let $\delta:=\delta_{n,T}$,
and $I_i:=\langle t_i,t_{i+1}]$. Then Cauchy-Schwarz inequality and isometry imply
\begin{equation}\begin{array}{lcl}
&&\E{\left|\frac{1}{T\sqrt{\delta}}\sum_{i=0}^{n-1}\int_{I_i}(C_{\vk}(X_t)\! -\! C_{\vk}(X_{t_i})) a(X_t)\,dt\right|}^2 =\\ 
&=&\frac{1}{T^2\delta}\E {\left|\int_0^T(\sum_{i=0}^{n-1}(C_{\vk}(X_t)\! -\! C_{\vk}(X_{t_i}))\1_{I_i}(t))\, a(X_t)\, dt\right|}^2\leq\\ 
&\leq&\frac{1}{T\delta}\E \int_0^T|\sum_{i=0}^{n-1}(C_{\vk}(X_t)\! -\! C_{\vk}(X_{t_i}))\1_{I_i}(t)|^2 a^2(X_t)\, dt =\\ 
&=&\frac{1}{T\delta}\E {\left|\int_0^T(\sum_{i=0}^{n-1}(C_{\vk}(X_t)\! -\! C_{\vk}(X_{t_i}))\1_{I_i}(t))\, a(X_t)\, dW_t\right|}^2=\\ 
&=&\E{\left|\frac{1}{\sqrt{T\delta}}\sum_{i=0}^{n-1}\int_{I_i}(C_{\vk}(X_t)\! -\! C_{\vk}(X_{t_i})) a(X_t)\,dW_t\right|}^2. 
\end{array}\label{iso}\end{equation}
Hence it is sufficient to prove that there exist constants $K_1>0$, $T_1\geq 0$, and $n_1$ such that
\begin{equation}\begin{array}{l}
\frac{1}{T\delta}\sum_{i=0}^{n-1}\!\!\int_{I_i}\!\!\E (|C_{\vk}(X_t)\! -\! C_{\vk}(X_{t_i})|^2 a^2(X_t))\,dt\leq
\! K_1\cdot K_\vk
\end{array}\label{novo}
\end{equation}
for $T>T_1$, and $n\geq n_1$ since the left hand side of (\ref{novo}) is equal to (\ref{iso}).
Similarly, to prove (\ref{ineq5}) and (\ref{ineq6})
it is sufficient to prove that there exist constants $K_2>0$, $T_0\geq T_1$, and $n_0\geq n_1$ such that for $T>T_0$, and $n\geq n_0$,
\begin{equation}\begin{array}{l}
\frac{1}{T\delta}\sum_{i=0}^{n-1}\!\!\int_{I_i}\!\!\E (|C_{\vk}(X_{t_i})|^2 {\left(\frac{b(X_t )}{b(X_{t_i})}-1\right)\!\!}^2 a^2(X_t))\,dt\leq
\! K_2\cdot K_\vk.
\end{array}\label{novo2}
\end{equation}

Let $j_1$,..., $j_d$  be nonnegative integers such that
$m:=j_1+\cdots +j_d\leq d+1$, and let $\vartheta\in\overline{\cal K}_0$ be fixed. Then function
$\tilde{f}:=\frac{\partial^m}{\partial\vartheta_1^{j_1}\cdots
\partial\vartheta_d^{j_d}}f(\cdot ,\vartheta)\in C^2(E)$ by (B2). If ${\cal A}\tilde{f}:=\tilde{f}'\mu_0 +\frac{1}{2}\tilde{f}''\nu^2$, then $|{\cal A}\tilde{f}|\leq g$, and $|\tilde{f}'\nu|\leq g$ by (B3).
Hence by applying It\^{o} formula, Jensen's inequality, and Lemma \ref{lema:4mom} it follows that
\[\begin{array}{lcl}
&&\E (\tilde{f}(X_t)-\tilde{f}(X_{t_i}))^4=\E (\int_{t_i}^t{\cal A}\tilde{f}(X_s)\,ds+\int_{t_i}^t(\tilde{f}'\nu)(X_s)\, dW_s)^4\leq\\
&\leq & 8(\E (\int_{t_i}^t|{\cal A}\tilde{f}|(X_s)\,ds)^4 +\E (\int_{t_i}^t(\tilde{f}'\nu)(X_s)\, dW_s)^4)\leq\\
&\leq & 8(\delta^3\, \E\int_{t_i}^t({\cal A}\tilde{f})^4(X_s)\,ds +3e^{3\delta} \E\int_{t_i}^t(\tilde{f}'\nu)^4(X_s)\,ds)\leq
24e^3 \E\int_{t_i}^t g^4(X_s)\, ds,
\end{array}\]
and $\E (a(X_t)-a(X_{t_i}))^4\leq 24e^3 \E\int_{t_i}^t \bar{a}^4(X_s)\, ds$ by an analogy, since we can assume that $\delta\leq 1$. Similarly,
\begin{equation}\begin{array}{lcl}
&&\E (a^2 (X_{t_i})(\tilde{f}(X_t)-\tilde{f}(X_{t_i}))^2)=\E [a^2 (X_{t_i})\E [(\tilde{f}(X_t)-\tilde{f}(X_{t_i}))^2|{\cal F}_{t_i}^0]]\leq\\ 
&\leq & 2\E [a^2 (X_{t_i})\E [\delta\,\int_{t_i}^t({\cal A}\tilde{f})^2(X_s)\,ds +(\int_{t_i}^t(\tilde{f}'\nu)(X_s)\, dW_s)^2|{\cal F}_{t_i}^0]]\leq\\ 
&\leq & 4\,\E [a^2 (X_{t_i})\int_{t_i}^t g^2(X_s)\, ds]\leq 2(t-t_i)\E a^4(X_{t_i})+2\E\int_{t_i}^t g^4(X_s)\, ds.
\end{array}\label{tilde_i}
\end{equation}
Hence
\begin{equation}\begin{array}{lcl}
&&\E((\tilde{f}(X_t)-\tilde{f}(X_{t_i}))^2 a^2(X_t))\leq\\
&\leq & \E (\tilde{f}(X_t)\!-\!\!\tilde{f}(X_{t_i}))^4+\E (a(X_t)\!-\! a(X_{t_i}))^4 +
2\E (a^2 (X_{t_i})(\tilde{f}(X_t)\!-\!\!\tilde{f}(X_{t_i}))^2)\leq\\
&\leq &25e^3\E\int_{t_i}^t g^4(X_s)\, ds + 24e^3\E\int_{t_i}^t\bar{a}^4(X_s)\, ds+4(t-t_i)\E a^4(X_{t_i}).
\end{array}\label{tilde}
\end{equation}
Now, let ${\vk}^{\vj}=k_1^{j_1}\cdots k_d^{j_d}$. Then
\begin{equation}\begin{array}{lcl}
&&{|\vk|}^{\vj\, 2}
\E (|C_{\vk}(X_t)\! -\! C_{\vk}(X_{t_i})|^2 a^2(X_t))=\E (|C_{\vk}^{(\vj)}(X_t)\! -\! C_{\vk}^{(\vj)}(X_{t_i})|^2 a^2(X_t))\leq\\ 
&\leq& \frac{1}{(2\pi)^d}\!\!\int_{\overline{\cal K}_0}\!\!\E{\left(\!\!{\left(\!\frac{\partial^m f}{\partial\vartheta_1^{j_1}\cdots
\partial\vartheta_d^{j_d}}(X_t,\vartheta)-\frac{\partial^m f}{\partial\vartheta_1^{j_1}\cdots
\partial\vartheta_d^{j_d}}(X_{t_i},\vartheta)\!\right)\!\!}^2\! a^2(X_t)\!\!\right)}\! d\vartheta\leq\\ 
&\leq &25e^3\E\int_{t_i}^t g^4(X_s)\, ds + 24e^3\E\int_{t_i}^t \bar{a}^4(X_s)\, ds+4(t-t_i)\E a^4(X_{t_i})
\end{array}\label{tilde2}
\end{equation}
by the definition of Fourier's coefficients, Jensen's inequality, Fubini's theorem, and (\ref{tilde}). Hence
\[\begin{array}{lcl}
&&(1+|k_1|+\cdots +|k_d|)^{2(d+1)} \E (|C_{\vk}(X_t)\! -\! C_{\vk}(X_{t_i})|^2 a^2(X_t))\leq\\
&\leq & (d+1)^{d+1}(1+|k_1|^2+\cdots +|k_d|^2)^{d+1} \E (|C_{\vk}(X_t)\! -\! C_{\vk}(X_{t_i})|^2 a^2(X_t))=\\
&=& (d+1)^{d+1}\!
\sum_{j_0+\cdots +j_d= d+1}\!\!\frac{(d+1)!}{j_0! j_1!\cdots j_d!}{|\vk|}^{\vj\, 2}\E (|C_{\vk}(X_t)\! -\! C_{\vk}(X_{t_i})|^2 a^2(X_t))\leq\\
&\leq & (d+1)^{2(d+1)}(25e^3\E\!\int_{t_i}^t\! g^4(X_s)\, ds + 24e^3\E\!\int_{t_i}^t\! \bar{a}^4(X_s)\, ds+4(t-t_i)\E a^4(X_{t_i}))
\end{array}\]
implying that 
\[\begin{array}{lcl}
&&\E (|C_{\vk}(X_t)\! -\! C_{\vk}(X_{t_i})|^2 a^2(X_t))\leq\\
&\leq & K_\vk^2 (25e^3\E\!\int_{t_i}^t\! g^4(X_s)\, ds + 24e^3\E\!\int_{t_i}^t\! \bar{a}^4(X_s)\, ds+4(t-t_i)\E a^4(X_{t_i})).
\end{array}\]
Finally, if $K'=5e^{3/2}$ then
\[\begin{array}{lcl}
&&\frac{1}{T\delta}\sum_{i=0}^{n-1}\!\!\int_{I_i}\!\!\E (|C_{\vk}(X_t)\! -\! C_{\vk}(X_{t_i})|^2 a^2(X_t))\,dt\leq \\
&\leq &K_\vk^2 {K'}^2\E(\frac{1}{T\delta}\sum_{i=0}^{n-1}\!\!\int_{I_i}\!\!\int_{t_i}^t\! (g^4+\bar{a}^4)(X_s)\, ds\, dt +
\frac{1}{T}\sum_{i=0}^{n-1} a^4(X_{t_i})\Delta_i t)\leq\\
&\leq & K_\vk^2 {K'}^2(\frac{1}{T}\E \int_0^T\!\!\!g^4(X_t)\,dt +\frac{1}{T}\E\!\!\int_0^T\!(\bar{a}^4(X_t)\,dt + 
\sum_{i=0}^{n-1} a^4(X_{t_i})\Delta_i t)).
\end{array}\]
Assumptions (B1-3) imply that there exist $T_1\geq 0$ and $n_1\in\N$ such that the expression in the parentheses on the right hand side of the above inequality is bounded by a constant ${K''}^2=C_g+C_a>0$ for $T>T_1$ and $n\geq n_1$. Hence $K_1:=K' K''>0$ in (\ref{novo}) and the statements (\ref{ineq3}-\ref{ineq4}) are proved. To prove (\ref{novo2}) and
hence statements (\ref{ineq5}-\ref{ineq6}) notice that
\[\begin{array}{l}
\frac{b(X_t)}{b(X_{t_0})}=\exp\left(\sqrt{\sigma}\int_{t_0}^t b'(X_s)\, dW_s+
\int_{t_0}^t\left(\frac{\mu_0 b'}{b} +\frac{\sigma}{2} (b''b-b'^2)\right)(X_s)\,ds\right)
\end{array}\]
from the proof of Lemma \ref{lema1}. It follows that
\begin{equation}\begin{array}{l}
\frac{b(X_t)}{b(X_{t_0})}-1=\int_{t_0}^t \frac{b(X_s)}{b(X_{t_0})}{\left(\sqrt{\sigma}b'(X_s)\, dW_s+
 (\frac{\mu_0 b'}{b} +\frac{\sigma}{2} b''b)(X_s)\,ds\right)\!}\label{itoexp}
\end{array}\end{equation}
by It\^{o} formula applied on the exponential function. Now, in the same way as equations (\ref{tilde_i}) and (\ref{tilde}) have been derived we obtain the following: for subdivisions such that $\delta_{n,T}\leq h_0$, and $\tilde{f}:=\frac{\partial^m}{\partial\vartheta_1^{j_1}\cdots
\partial\vartheta_d^{j_d}}f(\cdot ,\vartheta)$, $\vartheta\in\overline{\cal K}_0$, and some constants $C', C''>0$,
\[\begin{array}{lcl}
&&\E (a^2 (X_{t_i}) \tilde{f}^2 (X_{t_i}){\left(\frac{b(X_t)}{b(X_{t_i})}-1\right)\!\!}^2)\leq \E
(a^2 (X_{t_i}) g^2 (X_{t_i}){\left(\frac{b(X_t)}{b(X_{t_i})}-1\right)\!\!}^2)\leq \\
&\leq & C'( 2(t-t_i)\E (a^4(X_{t_i})g^4(X_{t_i}))+2\E\int_{t_i}^t{\left(\frac{b(X_s)}{b(X_{t_i})}\right)\!\!}^4 r^4(X_s)\, ds)\leq\\
&\leq & C'(2(t-t_i)\E (a^4(X_{t_i})g^4(X_{t_i}))+(t-t_i)\E\, c_0 (X_{t_i})+\E\int_{t_i}^t r^8(X_s)\, ds),
\end{array}\]
and,
\[\begin{array}{lcl}
&&\E (\tilde{f}^2 (X_{t_i}){\left(\frac{b(X_t)}{b(X_{t_i})}-1\right)\!\!}^2 a^2 (X_{t}))\leq\E (g^2 (X_{t_i}){\left(\frac{b(X_t)}{b(X_{t_i})}-1\right)\!\!}^2 a^2 (X_{t}))\leq\\
&\leq &\E (g^4 (X_{t_i})(a(X_{t})-a(X_{t_i}))^4)+\E{\left(\frac{b(X_t)}{b(X_{t_i})}-1\right)\!\!}^4 +\\
&&+2\, \E (a^2 (X_{t_i}) g^2 (X_{t_i}){\left(\frac{b(X_t)}{b(X_{t_i})}-1\right)\!\!}^2)\leq C''(\,
\E (g^4(X_{t_i})\int_{t_i}^t\!\bar{a}^4(X_s)\, ds)+\\
& &+(t-t_i)\E (a^4(X_{t_i})g^4(X_{t_i}))+(t-t_i)\E\, c_0 (X_{t_i})+\E\int_{t_i}^t r^8(X_s)\, ds)
\end{array}\]
by Lemma \ref{lema1}, (B1) and (B4). Hence there exist $K_2>0$, $T_0\geq T_1$, and $n_0\geq n_1$ such that for all
$T>T_0$ and $n\geq n_0$, (\ref{novo2}) follows in the same way as (\ref{novo}) has been followed from
(\ref{tilde_i}) and (\ref{tilde}) by using (B1-4).
\endproof\vspace{5pt}

{\em Proof of Lemma \ref{lema5}}. By Lemma \ref{lema4},
\[\begin{array}{l}
\sum_{\vk\in\Z^d}|C_{\vk}(x)|\leq  g(x)(d+1)^{d+1}
\sum_{\vk\in\Z^d}\frac{1}{(1+|k_1|+\cdots +|k_d|)^{d+1}},
\end{array}\]
and
\[\begin{array}{lcl}
&&\sum_{\vk\in\Z^d}\frac{1}{(1+|k_1|+\cdots +|k_d|)^{d+1}}=\\
&\leq &1+\sum_{r=1}^d{d\choose r}\frac{2^r}{r^{d+1}}\sum_{k_1=1}^\infty\cdots\sum_{k_r=1}^\infty\left(\frac{r}{k_1+\cdots +k_r}\right)^{d+1}\leq\\
&\leq& 1+\sum_{r=1}^d{d\choose r}\frac{2^r}{r^{d+1}}\sum_{k_1=1}^\infty\cdots\sum_{k_r=1}^\infty\frac{1}{
k_1^{\frac{d+1}{r}}\cdots k_r^{\frac{d+1}{r}}}=\\
&=&1+\sum_{r=1}^d{d\choose r}\frac{2^r}{r^{d+1}}\left(\sum_{k=1}^\infty\frac{1}{k^{\frac{d+1}{r}}}\right)^r\!\!.
\end{array}\]
Since $\sum_{k}(1/k^{\frac{d+1}{r}})<+\infty$ for all $r\leq d$, it follows that $K<+\infty$. Moreover, for any $N$, $\vartheta\in\overline{\cal K}_0$, and $x\in E$,
$
 |S_N(x,\vartheta )-f(x,\vartheta )|\leq
\sum_{|\vk|>N}|C_{\vk}(x)|\leq K g(x),
$
implying the statements of the lemma. \endproof\vspace{5pt}

{\em Proof of Lemma \ref{lema3}}. Let $x_0\in E$ be fixed and $F(x,\vartheta ):=\int_{x_0}^x\frac{f(y,\vartheta )}{\nu (y)} a(y)\, dy$. Then $F$
is a continuous function on $E\times\Theta$. By It\^{o} formula applied on $F$,
\[\begin{array}{lcl}
&&\int_0^T f(X_t,\vartheta )a (X_t)\, dW_t=\\
&=& F(X_T,\vartheta)\!-\!F(X_0,\vartheta)-\!\!\int_0^T\!\!\!\left(f(\cdot ,\vartheta)\!\!
\left(\frac{\mu_0}{\nu}-\frac{1}{2}\nu'\right)\! a+\frac{1}{2}(f (\cdot,\vartheta ) a)'\nu\right) (X_t)\,dt,
\end{array}\]
which is a continuous function on $\Theta$.\endproof\vspace{5pt}

{\em Proof of Lemma \ref{difkoef}}. Let $I_i:=\langle t_i ,t_{i+1}]$. Notice that
\[\begin{array}{cl}
&\frac{1}{T}\left|\sum_{i=0}^{n-1}\frac{1}{\Delta_i t}\left(\left(\int_{t_i}^{t_{i+1}}\frac{b(X_t)}{b(X_{t_i})}\, dW_t\right)^2 -(\Delta_i W)^2\right)\right|\leq\\
\leq & \frac{1}{T}\sum_{i=0}^{n-1}\frac{1}{\Delta_i t}(\int_{I_i}(\frac{b(X_t)}{b(X_{t_i})}-1) dW_t)^2 +\frac{2}{T}|\sum_{i=0}^{n-1}\frac{1}{\Delta_i t}\Delta_iW\int_{I_i}(\frac{b(X_t)}{b(X_{t_i})}-1) dW_t|.
\end{array}\]
Since
\[\begin{array}{l}
 \frac{1}{T}\E\sum_{i=0}^{n-1}\frac{1}{\Delta_i t}(\int_{I_i}(\frac{b(X_t)}{b(X_{t_i})}-1) dW_t)^2=
 \frac{1}{T}\sum_{i=0}^{n-1}\frac{1}{\Delta_i t}\int_{I_i}\E(\frac{b(X_t)}{b(X_{t_i})}-1)^2 dt
\end{array}\]
by the isometry, it follows that this expression is bounded by a constant for all $T>T_1$ and $n\geq n_1$ and
some $T_1\geq 0$ and $n_1$ in the same way as in the proof of Lemma \ref{lema6} since (B4) holds. It remains to prove the same for the second expression from the right hand side of the above inequality. By applying Ito formula
and (\ref{itoexp}) the following holds:
\[\begin{array}{cl}
&\Delta_iW\int_{I_i}(\frac{b(X_t)}{b(X_{t_i})}-1) dW_t =\\
=& \int_{I_i}\!(\int_{t_i}^t\!(\frac{b(X_s)}{b(X_{t_i})}-1) dW_s +(W_t -W_{t_i})(\frac{b(X_t)}{b(X_{t_i})}-1))\, dW_t
+\int_{I_i}\!(\frac{b(X_t)}{b(X_{t_i})}-1)\, dt=\\
=& \int_{I_i}\!\left(\int_{t_i}^t\!(\frac{b(X_s)}{b(X_{t_i})}-1) dW_s +(W_t -W_{t_i})(\frac{b(X_t)}{b(X_{t_i})}-1)+
\sqrt{\sigma}\Delta_i t\frac{b(X_t)}{b(X_{t_i})}b'(X_t)+\right.\\
& \left.+\sqrt{\sigma}(t-t_i)\int_{t_i}^t\!\frac{b(X_s)}{b(X_{t_i})}b'(X_s)
\, dW_s\right)\, dW_t+\int_{I_i}\!\int_{t_i}^t\!\frac{b(X_s)}{b(X_{t_i})}v(X_s)\, ds\, dt
\end{array}\]
where $v:=(\mu (\cdot ,\vartheta)/b)b' +(\sigma/2)b b''$. Then by applying the isometry and Cauchy inequality,
and by assuming that $T\geq 1$,
\[\begin{array}{cl}
&\E(\frac{1}{T}|\sum_{i=0}^{n-1}\frac{1}{\Delta_it}\Delta_iW\int_{I_i}(\frac{b(X_t)}{b(X_{t_i})}-1) dW_t|)^2 \leq\\
\leq & \frac{1}{T}\sum_{i=0}^{n-1}\!\left(\!2\frac{1}{\Delta_it^2}\!\int_{I_i}\!\int_{t_i}^t\!\E(\frac{b(X_s)}{b(X_{t_i})}-1)^2 ds\,dt+ \frac{1}{\Delta_it^2}\!\int_{I_i}\!\E(\frac{b(X_t)}{b(X_{t_i})}-1)^4dt+2\Delta_it +\right.\\
& +\sigma\int_{I_i}\E(\frac{b(X_t)}{b(X_{t_i})})^2r^2(X_t)dt+
\sigma\int_{I_i}\!\int_{t_i}^t\E(\frac{b(X_s)}{b(X_{t_i})})^2r^2(X_s)
\, ds\, dt+\\
&\left.+(1+\sigma/2)\frac{1}{\Delta_it}\int_{I_i}\!\int_{t_i}^t\E(\frac{b(X_s)}{b(X_{t_i})})^2r^2(X_s)\, ds\, dt\right)
\end{array}\]
since $|v|\leq (1+\sigma/2)r$ and $|b'|\leq r$ for function $r$ from (B4). For all terms on the right hand side of
the above inequality we can prove boundedness in the same way as in the proof of Lemma \ref{lema6} by using (B4), except for the following one for which we have to use the additional assumptions of the lemma to obtain the boundedness. First by using (\ref{itoexp}), then Lemma \ref{lema:4mom}, and Ito formula we obtain the following:
\[\begin{array}{cl}
&\frac{1}{T}\sum_{i=0}^{n-1}\!\frac{1}{\Delta_it^2}\!\int_{I_i}\E(\frac{b(X_t)}{b(X_{t_i})}-1)^4dt\leq\\
\leq & K'(1+\frac{1}{T}\E(\int_0^T (r^8+r^{16}+(b^2b''')^8)(X_t)\,dt+\sum_{i=0}^{n-1}(c_0+r^4)(X_{t_i})\Delta_it)
\end{array}\]
for some constant $K'>0$. Now, the statement of the lemma follows.
\endproof\vspace{5pt}

{\em Proof of Lemma \ref{lm:25}}. By applying (\ref{e1}) it follows that:
\[\begin{array}{lclr}
&&\!\! \sum_{i=0}^{n-1}\frac{(\Delta_i X -\mu (X_{t_i},\vartheta )\Delta_it )^2}{
\sigma b^2 (X_{t_i})\Delta_it}-\sum_{i=0}^{n-1}\frac{(\Delta_i W)^2}{\Delta_it}= & \\
&=&\!\!\sum_{i=0}^{n-1}\frac{1}{
\Delta_it}\left(\int_{t_i}^{t_{i+1}}\!\!\frac{\mu (X_t,\vartheta)\! -\!
\mu(X_{t_i},\vartheta)}{b(X_{t_i})}\,dt\right)^2- & (\mbox{E1})\\ 
&&\!\!-2\sum_{i=0}^{n-1}\!\!\sqrt{\sigma}\int_{t_i}^{t_{i+1}}\!\!\frac{b(X_t)}{b(X_{t_i})}\, dW_t\cdot\frac{1}{
\Delta_it}\!\!\int_{t_i}^{t_{i+1}}\!\!\frac{\mu (X_t,\vartheta)\!-\!
\mu(X_{t_i},\vartheta)}{b(X_{t_i})}\,dt+ & (\mbox{E2})\\ 
&&\!\!+\sigma\sum_{i=0}^{n-1}\frac{1}{
\Delta_it}\left((\int_{t_i}^{t_{i+1}}\!\!\frac{b(X_t)}{b(X_{t_i})}\, dW_t)^2-
(\Delta_iW)^2\right). & (\mbox{E3}) 
\end{array}\]
First we will prove that
the expression from the left hand side of the above equation is bounded in $L^1$-norm by a constant for all $n\geq n_0$ (for some $n_0$) in case when
all functions $\mu (\cdot,\vartheta )$, $b$ and their appropriate partial derivatives are bounded on $E$, and then the statement of the lemma will follow by using local compactness of $E$ and Markov's inequality just in the same way as in the proof of Corollary \ref{tm:23}. Let $f:=\mu(\cdot,\vartheta )/b$ and $I_i:=\langle t_i,t_{i+1}]$, and let $n$ be such that $\delta_{n,T}\leq 1$. Then the expectation of $(E1)$ is dominated by
\[\begin{array}{lcl}
\!\!& &\E\! \sum_{i=0}^{n-1}\!\!\frac{1}{\Delta_it}\!\left(\!\int_{I_i}\!\!\frac{\mu (X_t,\vartheta)-
\mu(X_{t_i},\vartheta)}{b(X_{t_i})}\,dt\right)^2\!\!\!\leq 
\!\E \!\sum_{i=0}^{n-1}\!\!\frac{1}{(\Delta_it)^2}\!{\left(\int_{I_i}\!\!\frac{\mu (X_t,\vartheta)-
\mu(X_{t_i},\vartheta)}{b(X_{t_i})}\,dt\right)\!}^2\!\!\leq \\
\!\!&\leq& 2T \left(\frac{1}{T}\sum_{i=0}^{n-1}\frac{1}{\Delta_it}\int_{I_i}
\E (f(X_t)-
f(X_{t_i}))^2\,dt +\right.\\
\!\!& & \left.+\frac{1}{T}\sum_{i=0}^{n-1}\frac{1}{\Delta_it}\int_{I_i}
\E f^2(X_t)(\frac{b(X_t)}{b(X_{t_i})}-1)^2\,dt\right)\!\leq T C'.
\end{array}\]
The existence of a constant $C'>0$  follows in the same way as in the proof of Lemma \ref{lema6} since (B1-4) hold
for bounded functions by Remark \ref{rem3}. $L^1$-norm of $(E2)$ is dominated by
\[\begin{array}{cl}
&\E\!\sum_{i=0}^{n-1}\left(\int_{t_i}^{t_{i+1}}\!\!\frac{b(X_t)}{b(X_{t_i})}\, dW_t\right)^2+
\E \!\sum_{i=0}^{n-1}\!\!\frac{1}{(\Delta_it)^2}\!{\left(\int_{I_i}\!\!\frac{\mu (X_t,\vartheta)-
\mu(X_{t_i},\vartheta)}{b(X_{t_i})}\,dt\right)\!}^2\leq\\
\leq &
TK(1+\frac{1}{n}\sum_{i=0}^{n-1}c_0(X_{t_i}))+TC'\end{array}\]
for some constant $K>0$ by the isometry and Lemma \ref{lema1}. Now, boundedness of $L^1$-norm of $(E2)$ follows from $(B4)$.
$L^1$-norm of $(E3)$ is bounded by Lemma \ref{difkoef} and Remark \ref{rem4}.
\endproof\vspace{5pt}

{\em Proof of Lemma \ref{lema1b}}. Let $f:=\mu/b:E\times\overline{\Theta}\rightarrow\R$, and let
$f_0:=\mu (\cdot,\vartheta_0)/b$. 
For  nonnegative integers
$j_1$,..., $j_d$ such that $j_1+\cdots +j_d=3$, let
$\tilde{f}:=\frac{\partial^3}{\partial\vartheta_1^{j_1}\cdots\partial\vartheta_d^{j_d}}f$, and
$\hat{f}:=\frac{\partial^3}{\partial\vartheta_1^{j_1}\cdots\partial\vartheta_d^{j_d}}(f^2)$.  Then for $T>0$,
and $\vartheta\in\overline{\Theta}$,
\[\begin{array}{l}
\frac{\partial^3}{\partial\vartheta_1^{j_1}\cdots\partial\vartheta_d^{j_d}}\ell_T (\vartheta )=
\int_0^T (\tilde{f}(\cdot,\vartheta )f_0-\frac{1}{2}\hat{f}(\cdot,\vartheta ))(X_t)\, dt+\sqrt{\sigma}\int_0^T
\tilde{f}(X_t,\vartheta )\, dW_t
\end{array}\]
by (\ref{Dell}) and (\ref{e1}). Since {\sc (H2b-3b)} hold,  from the proof of Corollary \ref{tm:cor2} it follows that
\[\begin{array}{l}
\sup_{\vartheta\in\overline{\Theta}}\frac{1}{T}|\int_0^T (\tilde{f}(\cdot,\vartheta )f_0-\frac{1}{2}\hat{f}(\cdot,\vartheta ))(X_t)\, dt|\leq C\frac{1}{T}\int_0^T g_{00}^2(X_t)\, dt,
\end{array}\]
 where $C:=1+7\cdot 2^3$. The right hand side of the above inequality $\Pb_{\theta_0}$-a.s. converge to a finite nonrandom limit $L_0=L_0(\vartheta_0)$ by the ergodic property of $X$. Hence on an $\Pb_{\theta_0}$-a.s. event there exists $T_0'\geq 0$
 such that for all $T >T_0'$,
 \[\begin{array}{l}
 \frac{1}{T}\int_0^T g_{00}^2(X_t)\, dt\leq |\frac{1}{T}\int_0^T g_{00}^2(X_t)\, dt-L_0|+L_0\leq 1+L_0.
 \end{array}\]
 Let us suppose that $b>0$. The case when $b<0$ can be analyzing in the same way. By applying It\^{o} formula twice, first on function $x\mapsto\int_{x_0}^x \frac{\tilde{f}(y,\vartheta)}{b(y)}\, dy$, and then
on $x\mapsto\int_{x_0}^x \frac{g_0(y)}{b(y)}\, dy$, we get the following
\[\begin{array}{cl}
&\frac{\sqrt{\sigma}}{T}|\!\int_0^T\!\!
\tilde{f}(X_t,\vartheta )\, dW_t|=|\frac{\sqrt{\sigma}}{T}\!\int_{x_0}^{X_T}\!\frac{\tilde{f}(y,\vartheta)}{b(y)}\, dy\!-\!\frac{1}{T}
\!\int_0^T\!(\tilde{f}f_0\!-\!\frac{\sigma}{2}(\tilde{f}b'\!-\!\tilde{f}'b))(X_t)\, dt|\leq\\
\leq\!&\frac{\sqrt{\sigma}}{T}|\!\int_{x_0}^{X_T}\!\frac{g_0(y)}{b(y)}\, dy|+\frac{1}{T}\int_0^T
(g_0^2+\frac{\sigma}{2}(g_1+2g_0|b'|))(X_t)\, dt\leq\\
\leq\!&\frac{\sqrt{\sigma}}{T}|\!\int_{0}^{T}\!g_0(X_t)\, dW_t|\!+\!\frac{1}{T}
\!\int_0^T\!(2g_0^2\!+\!\frac{\sigma}{2}(g_0\!+\!|g_0'b|+g_1+2g_0|b'|))(X_t)\, dt.
\end{array}\]
Since the right hand side of the above inequality $\Pb_{\theta_0}$-a.s. converge to a finite nonrandom limit (by ergodic property and the law of large numbers for continuous martingales since {\sc (H2b-3b)} hold), and since it is also
an upper bond for the left hand side uniformly for all $\vartheta\in\overline{\Theta}$, and all partial derivatives of the third order, there exists a constant $C''>0$ such that $\Pb_{\theta_0}$-a.s. there exists $T_0''\geq T_0'$ such that for all $T>T_0''$,
$\sup_{\vartheta\in\overline{\Theta}}\frac{1}{T}|D^3\ell_T (\vartheta )|_\infty\leq C''$. From the definition of operator norm, for the same $T$,
\[\begin{array}{l}
\sup_{\vartheta\in\overline{\Theta}}\frac{1}{T}|D^3\ell_T (\vartheta )|\leq \sup_{\vartheta\in\overline{\Theta}}\frac{1}{T}d^{3/2}|D^3\ell_T (\vartheta )|_\infty\leq d^{3/2}C''=:C_2.
\end{array}\]
By the same arguments we can prove that there exist constants $C_0>0$, $C_1>0$ such that $\Pb_{\theta_0}$-a.s.
exists $T_0\geq 0$ such that $T_0\geq T_0''$, and
 for all $T>T_0$, $\sup_{\vartheta\in\overline{\Theta}}\frac{1}{T}|D^{r+1}\ell_T (\vartheta )|\leq C_r$ for $r=0,1$. Finally, the statements of the lemma follow from the mean value theorem and Taylor expansion
(\ref{eltaylor}) from the proof of Theorem \ref{tm:cont}, where $\hat{\vartheta}_T$ and $\vartheta_0$ are replaced with
 $\vartheta_1$ and $\vartheta_2$ respectively.\endproof\vspace{5pt}

{\em Proof of Lemma \ref{lema2b}}. Let $C_0>0$ and $\Omega_0$ be an event, both from Lemma \ref{lema1b} such that $\Pb_{\theta_0}(\Omega_0)=1$ and on $\Omega_0$ for all $T\geq T_0$, and all
$\vartheta_1,\vartheta_2\in\overline{\Omega}$, $|\ell_T(\vartheta_1)-\ell_T(\vartheta_2)|\leq TC_0|\vartheta_1-\vartheta_2|$. Let $K_0>0$ be Lipschitz constant of function $\ell_{\vartheta_0}$, and let
$\varepsilon >0$ be an arbitrary number. Let $\delta:=\varepsilon/(2(C_0+K_0))$. Since $\{ K(\vartheta,\delta ):\vartheta\in\Theta\}$ is an open cover of compact $\overline{\Theta}$, there exists a finite subcover
$\{K(\vartheta_i,\delta):i=1,\ldots, K_\varepsilon\}$. Let $\Omega_1$ be an $\Pb_{\theta_0}$-a.s. event such
that on this event there exists $T_\varepsilon\geq T_0$ such that for all $T\geq T_\varepsilon$,
and $1\leq j\leq K_\varepsilon$,
$|\frac{1}{T}\ell_T(\vartheta_j)\!-\!\ell_{\vartheta_0}(\vartheta_j)|<\varepsilon/(2K_\varepsilon)$. Then
on $\Omega_0\cap\Omega_1$ for all $\vartheta\in\Theta$ there exists $i=i(\vartheta)\leq K_\varepsilon$ such that
$\vartheta\in K(\vartheta_i,\delta)$, and
\[\begin{array}{cl}
&\!|\frac{1}{T}\ell_T(\vartheta)\!-\!\ell_{\vartheta_0}(\vartheta)|\!=\!|\frac{1}{T}\ell_T(\vartheta)\!-\!
\frac{1}{T}\ell_T(\vartheta_i)\!+\!\frac{1}{T}\ell_T(\vartheta_i)\!-\!\ell_{\vartheta_0}\!(\vartheta_i)\!+
\!\ell_{\vartheta_0}\!(\vartheta_i)\!-\!\ell_{\vartheta_0}\!(\vartheta)|\!\leq\\
\leq\! & C_0 |\vartheta-\vartheta_i|\!+\!\sum_{j=1}^{K_\varepsilon}|\frac{1}{T}\ell_T(\vartheta_j)\!-\!\ell_{\vartheta_0}(\vartheta_j)|
\!+\! K_0|\vartheta-\vartheta_i|<\\
<\! & (C_0+K_0)\delta + K_\varepsilon\frac{\varepsilon}{2K_\varepsilon}= \varepsilon.
\end{array}\]
Hence $\sup_{\vartheta\in\overline{\Theta}}|\frac{1}{T}\ell_T(\vartheta)\!-\!\ell_{\vartheta_0}(\vartheta)|<\varepsilon$ which proves the lemma.\endproof
}

\end{document}